\documentclass[12pt,a4paper,leqno]{article}
\usepackage {amsmath, latexsym, amscd, amssymb}
\usepackage {graphicx}
\usepackage{mathabx}
\usepackage{lineno}
\def\tvi #1#2{\vrule height #1pt depth #2pt width 0mm}

\newcommand{\N}{\mathbb{N}}

\newcommand{\C}{\mathbb{C}} 
\newcommand{\R}{\mathbb{R}}
\newcommand{\Z}{\mathbb{Z}}

\newcommand{\sB}{\mathcal{B}}
\newcommand{\sC}{\mathcal{C}}

\newcommand{\sH}{\mathcal{H}}

\newcommand{\sM}{\mathcal{M}}

\newcommand{\sN}{\mathcal{N}}
\newcommand{\sO}{\mathcal{O}}
\newcommand{\stO}{\mathcal{\widetilde O}}
\newcommand{\stH}{\mathcal{\widetilde H}}

\newcommand{\sR}{\mathcal{R}}
\newcommand{\sT}{\mathcal{T}}

\newcommand{\ugt}{{\theta^\star}}
\newcommand{\ugo}{{\omega^\star}}

\newcommand{\cf}{{\it cf.} }
\newcommand{\ie}{{\it i.e.}, }

\newcommand{\proof}{ \noindent { \em {\sc Proof.} }}
\newcommand{\qed}{ $\hfill {\Box}$\break}

\newcommand{\diag}{{\rm diag} }

\newcommand{\ga }{\alpha }

\newcommand{\gc }{\gamma }
\newcommand{\gd }{\delta }

\newcommand{\gve }{\varepsilon }
\newcommand{\gvf }{\varphi }

\newcommand{\gk }{\kappa }
\newcommand{\gz }{\zeta }
\newcommand{\gl}{\lambda}

\newcommand{\go }{\omega }

\newcommand{\gt}{\theta}

\newcommand{\gC }{\Gamma }
\newcommand{\gD }{\Delta }
\newcommand{\gO }{\Omega }
\newcommand{\gL}{\Lambda}

\newcommand{\gS}{\Sigma}

\newcommand{\hB}{\widehat B}

\newcommand{\hy}{\widehat y}
\newcommand{\hF}{\widehat F}
\newcommand{\hH}{\widehat H}
\newcommand{\hf}{\widehat f}
\newcommand{\hg}{\widehat g}
\newcommand{\hh}{\widehat h}

\def\hfpq #1{\widehat f^{\scriptscriptstyle{\bullet};q}}
\def\hBjp #1{\widehat B^{#1;\scriptscriptstyle{\bullet}}}
\def\vBjp #1{\vert\widehat B\vert^{#1;\scriptscriptstyle{\bullet}}}
\def\gvfpq #1{\gvf^{\scriptscriptstyle{\bullet};#1}}
\def\phipq #1{\phi^{\scriptscriptstyle{\bullet};#1}}
\def\phijp #1{\gvf^{#1;\scriptscriptstyle{\bullet}}}

\def\tfjp #1{\widetilde f^{#1;\scriptscriptstyle{\bullet}}}
\def\tgjp #1{\widetilde g^{#1;\scriptscriptstyle{\bullet}}}
\def\tfjq #1#2{\widetilde f^{#1;#2}}
\def\tfpq #1{\widetilde f^{\scriptscriptstyle{\bullet};#1}}

\def\tvi #1#2{\vrule height #1pt depth #2pt width 0mm}

\newcommand{\tY}{\widetilde Y}

\newcommand{\tH}{\widetilde{\sl H}}

\newcommand{\tF}{\widetilde F}
\newcommand{\tf}{\widetilde f}
\newcommand{\tg}{\widetilde g}

\newcommand{\Cf}{\C\lbrack\lbrack x\rbrack\rbrack}
\newcommand{\Cc}{\C\lbrace x\rbrace}

\newcommand{\dD}{\stackrel{\scriptscriptstyle\bullet}{\Delta}}

\newcommand{\dDodp}{\stackrel{\scriptscriptstyle \, \bullet\, + }{\Delta_{\ugo}}}
\newcommand{\ddodp}{\stackrel{\scriptscriptstyle \, \bullet\, + }{\delta_{\ugo}}}

\newtheorem{lemma}{Lemma}[section]
\newtheorem{thm}[lemma]{Theorem}
\newtheorem{cor}[lemma]{Corollary}
\newtheorem{pro}[lemma]{Proposition}
\newtheorem{dfn}[lemma]{Definition}

\newtheorem{rem}[lemma]{Remark}

\usepackage{hyperref}

\title{{\bf 
Resurgence, Stokes phenomenon and alien derivatives for level-one linear differential systems. \\
}}
\author{Mich\`ele LODAY-RICHAUD \& Pascal REMY}

\begin{document}




 \baselineskip=1.3em
\maketitle
\begin{abstract}
A precise description of the singularities of the Borel transform of solutions of a level-one linear differential system is deduced from a proof of the summable-resurgence of the solutions by the perturbative method of J. \'Ecalle. Then we compare the meromorphic classification (Stokes phenomenon) from the  viewpoint of the Stokes cocycle and the viewpoint of alien derivatives. We make explicit the Stokes-Ramis matrices as functions of the connection constants in the Borel plane and we develop two examples. No assumption of genericity is made.
\end{abstract}
{\let\thefootnote\relax\footnotetext{{\bf Keywords}:
Linear differential systems, linear differential equations, Stokes phenomenon, summability, resurgence, alien derivatives.

{\bf \ AMS classification}: 34M03, 34M30, 34M35, 34M40.
}}
\tableofcontents


\section{Introduction}\label{intro}

All along the paper we are given an ordinary linear
differential system (in short, a differential system or a 
system) of dimension $n$  with analytic coefficients at 0 
in $\C$ and rank one\footnote{The  rank is the order of the pole $x=0$ minus 1  in the system written in ``solved form'' $\displaystyle \frac{dY}{dx}=\frac{1}{x^2}A(x)\,Y$.} 
\begin{equation}\label{A}
x^2\frac{dY}{dx}= A(x)Y \qquad A(x)\in M_n(\Cc), \ A(0)\neq 0
\end{equation}
together with a formal fundamental solution at 0
\begin{equation}\label{ffs}
\tY(x)=\tF(x) x^L e^{Q(1/x)}
\end{equation}
where 
$$\begin{array}{ll}
\bullet \ \tF(x)\in M_n\big(\C(( x))\big),& \text{ is an invertible formal meromorphic matrix,}\\
\noalign{\medskip}
\bullet \displaystyle\ L=\bigoplus_{j=1}^J (\gl_j I_{n_j}+J_{n_j})& 
\text{ where } I_{n_j} \text{ is the identity matrix of size } n_j\\
&  \text{ and } J_{n_j}=\begin{bmatrix}
 0&1&\cdots&0\\
 \vdots&\ddots&\ddots&\vdots\\
 \vdots&&\ddots&1\\
 0&\cdots&\cdots&0\\
\end{bmatrix}
\text{ is an irreducible}\\
\noalign{\smallskip}
& \text{  Jordan block of size } n_j \ (J_{n_j}=0 \text{ if } n_j=1),\\
\end{array}
$$
$$
\begin{array}{ll}
\bullet\displaystyle\  Q(1/x)=\bigoplus_{j=1}^J q_j(1/x) I_{n_j}
&\quad  \text{ where the } q_j's \text{ are polynomials}.\\
\end{array}
$$
In the very general rank one case the determining polynomials $q_j$ are of maximal degree equal to 1 with respect to $1/x$ but they could be polynomials in a fractional power of $1/x$. Our assumption of ``single level equal to 1''
implies that the polynomials $q_j$ be monomials of degree 1
in $1/x$, not all  equal to a same polynomial $q$, some of them being possibly zero. We denote  
\begin{equation}\label{nivpur1}
 Q(1/x)=\bigoplus_{j=1}^J \frac{-a_j}{x} I_{n_j}.
\end{equation}

The system
\begin{equation}\label{A0}
x^2\frac{dY}{dx}=A_0(x)\,Y
\end{equation}
with formal fundamental solution $\displaystyle \tY_0(x)=x^L\,e^{Q(1/x)}$  has analytic coefficients and is called {\it a normal form} of System (\ref{A}). The fundamental solution $\tY_0(x)$ is called {\it a normal solution}.
It provides all formal invariants of System (\ref{A}), \ie invariants under formal gauge transformations $Y\mapsto \widetilde \Phi(x) Y$ with $\widetilde \Phi\in GL_n\big(\C(( x))\big)$
\footnote{The formal classification over an extension $\C(( t))$ of the base field $\C(( x))$ by a ramification $\displaystyle x=t^p$ started with Poincar\'e and Fabry and later Turrittin  \cite{T55}. See also Malgrange \cite{DMR07}. Over the base field itself it is due to Balser, Jurkat and Lutz \cite{BJL79}. For a shorter proof and an improvement of the Formal Classification Theorem see \cite{L-R01}.}.
\smallskip

Note that, in the whole paper, we make no other assumption  than the assumption of a single level equal to one, case in which the basic theory of resurgence takes place. In particular, 
no assumption of genericity, such as distinct eigenvalues or 
diagonal monodromy $L$, is made and any kind of resonance is allowed.
\medskip

The paper deals with the analytical properties of the solutions of System (\ref{A}) underlying the meromorphic classification at 0, also referred to as {\it Stokes phenomenon}. There exists mainly two ``dual'' approaches to this phenomenon:
\begin{itemize}
\item The first one , related to the theory of summation, is fully developed in the plane of the initial variable $x$, which we refer to as the {\it Laplace plane}. 
Various methods 
\cite{Mal83, MR91, BV89, L-R94} produce a full set of invariants in the form of 
{\it Stokes matrices}. 
Are considered only those Stokes matrices providing the transition between the sums (in our case, Borel-Laplace or 1-sums)  of a same formal fundamental solution $\tY(x)$ on each side of its 
anti-Stokes (singular) directions. We call them {\it Stokes-Ramis matrices}.
\item  The second one is strongly related  
to the {\it theory of resurgence} \cite{Ec81} and produces invariants in the form of {\it alien derivations}. The alien derivatives of a series $\tf(x)$
have been  defined by J. \'Ecalle through an average of various analytic continuations of 
the Borel transform $\hf(\xi)$ of $\tf(x)$; they have been  mostly developed in the plane of the variable $\xi$ which we refer to as the {\it Borel plane}.  
\end{itemize}

Much have been already said on these questions but either in situations restricted by generic conditions or in very general ones. In particular, the theory of resurgence was developed by J. \'Ecalle in the very general framework of non-linear differential equations, difference equations and so on\dots where it proves to be very efficient. In this context it seemed to us to be worth to make explicit what is specific to the linear case, what has to be really taken under consideration or can be made more precise, and how the various viewpoints are connected.

\bigskip

The first Section is devoted to  proving that the formal gauge transformation $\tF(x)$ is summable-resurgent. We first sketch a proof on linear differential equations based on the Newton polygon and Ramis Index Theorems. We develop then a proof by perturbation and majorant solutions following  J. \'Ecalle's method \cite{Ec85}.

This latter proof allows us to display a precise description of the singularities in the Borel plane (Section 2).

In Section 3 we compare the classification by the  Stokes cocycle versus the classification by the alien derivations. Roughly speaking,  the first approach selects an element in a unipotent Lie group while the second one provides its ``tangent'' variant in the associated Lie algebra. Moreover, for theoretical reasons as well as for computational ones, we make explicit the Stokes matrices in terms of the  connection constants in the Borel plane.

\vspace{.5cm}



\section{Summable-resurgence}  
\subsection{Prepared system}\label{prepared}
Before to start the calculations we prepare the system as follows.

A gauge transformation of the form $Y \mapsto T(x)\,x^{-\gl_1}\, e^{a_1/x} Y$ 
where the transformation $T(x)$ has explicit computable 
polynomial entries in $x$ and $1/x$ allows to assume that the 
following conditions are satisfied:
\begin{equation}\label{F=I+}
\tF(x)=I_n+\sum_{m\geq 1}F_mx^m\in GL_n(\Cf) \text{ with initial 
condition } 
 \tF(0)=I_n,\\
\end{equation}
\vspace{-.9em}
\begin{equation} \label{Re}
0\leq {\rm Re} (\gl_j)<1 \ \text{ for } \ j=1,\dots,J,
\end{equation}
\vspace{-.9em}
\begin{equation}\label{a1-0}
a_1=\gl_1=0.
\end{equation}

Conditions (\ref{F=I+}) and (\ref{Re}) guaranty the unicity 
of $\tF(x)$. Condition (\ref{a1-0}) is for notational convenience. 
Still, the $a_j$'s are not supposed distinct.

Any of the $J$ column-blocks of $\tF(x)$ associated with the 
irreducible Jordan blocks of $L$ (matrix of exponents of 
formal monodromy) can be positioned at the first place by means 
of a permutation $P$ on the columns of $\tY(x)$.  If $\tY(x)$ is given
in the form (\ref{ffs}) so is the new formal fundamental solution
$\tY(x)P=\tF(x)P x^{P^{-1}LP}e^{P^{-1}{Q(1/x)}P}$. Thus, we can  
restrict our study to the first column-block of $\tF(x)$ that we
denote by $\tf(x)$. 

Gauge transformations and permutations of this kind will be referred to as {\it elementary transformations}.


\subsection{Some definitions}\label{knownresults}
For the convenience of the reader we recall some definitions about the notions of resurgence and summation adapted to linear level-one differential systems. It is worth to note that, due to the linearity, it is useless, at least for the moment,  to consider
convolution algebras and lattices of singularities since no
convolution of singularities may occur. 
For a more general framework  we refer to \cite{Ec85, Mal85, Sau05}. 
\bigskip

All along the article, given a matrix $M$ split into blocks fitting the structure of $L$, we denote
\begin{itemize}
\item $M^{j;\scriptscriptstyle{\bullet}}$ the $j^{th}$ row-block of $M$,
\item $M^{\scriptscriptstyle{\bullet};k}$ the $k^{th}$ 
column-block of $M$,
\item $M^{j;k}$ the $k^{th}$ column-block in the $j^{th}$ 
row-block of $M$,
\item $M^{j;(k)}$ the  $k^{th}$ column in the $j^{th}$ row-block of $M$,
\item $M^{(j,\ell);\scriptscriptstyle{\bullet}}$ the $\ell^{th}$ row in the 
$j^{th}$ row-block of $M$.
\end{itemize}

Let $\gO=\{a_j, \, j=1,\dots, J\}$ denote the set of {\it Stokes values} 
associated with  System (\ref{A}).
Theorem \ref{thmsum-res} below asserts, in particular, that all  
possible singularities of the Borel transform $\hf(\xi)$  of $\tf(x)$
belong to $\gO$.  The directions determined by the
elements of $\gO^*=\gO\setminus\{0\}$ 
from 0 are called  {\it anti-Stokes directions associated with}
$\tf(x)$.  Given a direction $\gt\in \R/2\pi\Z$ let $d_\gt$ denote the half line issuing from 0 with argument $\gt$. We denote 
\begin{itemize}
\item $\gO_\gt=\gO^*\cap d_\gt$ the set of non-zero Stokes values of System (\ref{A})  with argument $\gt$.
\end{itemize}
The anti-Stokes directions associated with the $k^{th}$ column-block of $\tF(x)$ are given by the non-zero elements of $\gO-a_k$ (to normalize the $k^{th}$ column-block one has to multiply by $\displaystyle e^{-a_k/x} $) and the {\it anti-Stokes directions of System} (\ref{A}) (\ie associated with the full matrix $\tF(x)$) are given by the non-zero elements~of 
\begin{itemize}
\item ${\bf \gO}=\{\, a_j-a_k\, ;\, 1\leq j,k\leq J\}.$
\end{itemize}
The elements of ${\bf \gO}$ are the Stokes values of the homological  system satisfied by $\tF(x)$ (\cf System (\ref{AA0}) below).\\
We denote
\begin{itemize}
\item ${\bf \gO}_\gt={\bf \gO}^*\cap d_\gt$ the set of non-zero Stokes values of System (\ref{AA0}) with argument $\gt$.
\end{itemize}
\bigskip

The adequate Riemann surface  on which the Borel
transform $\hf(\xi)$ of $\tf(x)$ lives is the surface $\sR_{\gO}$ 
 defined as below. 

\begin{dfn}\label{RS0} {\bf\protect\boldmath Riemann surface
    $\sR_{\gO}$.\protect\unboldmath } 
\begin{itemize}
\item {\rm The points of  $\sR_{\gO}$ are the homotopy
  classes with fixed extremities of paths $\gc$ issuing from $0$ and
  lying in $\C\setminus \gO$ (but the starting point $0$); in
  particular, no path
  except those that are homotopic to the constant path $0$ ends 
  at $0$. }
\item {\rm The complex structure of  $\sR_{\gO}$ is the pull-back of the
  usual complex structure of $\C\setminus \gO^*$ by the natural 
  projection $$\displaystyle   
\pi : \left\lbrace\begin{array}{ccl}
 \sR_{\gO}&\longrightarrow& \C\setminus \gO^* \\
 \lbrack \gc\rbrack&\mapsto&\text{end-point } \gc(1).\\
\end{array}\right.$$}
\end{itemize}
\end{dfn}

The difference between  $\sR_{\gO}$ and the universal cover of
$\C\setminus \gO$ lies in the fact that  $\sR_{\gO}$ has no branch point at $0$ in the first sheet. 
\smallskip

\begin{dfn} \label{res} {\bf Resurgence.}
\begin{itemize}
\item {\rm \underline{Resurgence in the Borel plane} (in Ecalle's
  language, {\it convolutive model}).\\ 
A {\it resurgent function} with singular support $\gO$  is any
function defined and analytic on all of $\sR_{\gO}$. }

\item {\rm\underline{Resurgence in the Laplace plane} (in Ecalle's
  language, {\it formal model}).\\ 
A  series $\tf(x)$ in powers of $x$ is said to be a {\it resurgent series} with
singular support $\gO$  when its Borel transform $\hf(\xi)$  is
a convergent series near $\xi=0$ which can be analytically continued to $\sR_{\gO}$ (in
short, its Borel transform is resurgent with
singular support $\gO$).
\\ 
{\rm Let $\widehat{\sR es}_{\gO}\tvi{15}{0}$  and $\widetilde{\sR es}_{\gO}\tvi{15}{0}$ denote the sets of resurgent functions and  of resurgent series with singular
support $\gO$ respectively. }}
\end{itemize}
\end{dfn}
A resurgent function of $\widehat{\sR es}_{\gO}\tvi{15}{0}$ is analytic at $0$ on the first sheet of $\sR_{\gO}$. To emphasize the special role played by 0 we sometimes denote the singular support by $\gO,0$
As a Fourier operator the Borel transform exchanges multiplication by an exponential and translation. We denote by $\sR_{\gO,\go} $ and $\widehat{\sR es}_{\gO,\go}$  the spaces $\sR_{\gO} $ and $\widehat{\sR es}_{\gO}$ translated by $\go$ so that, in particular, $\sR_{\gO} =\sR_{\gO,0} $ and $\widehat{\sR es}_{\gO}=\widehat{\sR es}_{\gO,0}$. If a series $\tf(x)$ belongs to $\widetilde{\sR es}_{\gO}$ the Borel transform of $\displaystyle \tf(x)e^{-\go/x}$ belongs to $\widehat{\sR es}_{\gO,\go}$.

\bigskip

While in non linear situations it is soon necessary to endow $\widehat{\sR  es}_{\gO}$  with a structure of a convolution algebra ---usually a quite
difficult task--- in the linear case under consideration  we are going to meet, at least  in this section,  only 
the convolution of elements of $\widehat{\sR es}_{\gO}$ with entire
functions over $\C$ growing at most exponentially at infinity.  Let
$\sO^{\leq 1}(\C)$ denote the convolution algebra of entire functions
on $\C$ with exponential growth at infinity.  The set
$\widehat{\sR es}_{\gO}$  has a natural structure of a $\sO^{\leq
  1}(\C)*$module (the star refers to the fact that multiplication by a
scalar is the convolution by an 
element of $\sO^{\leq 1}(\C)$). There  corresponds, on $\widetilde{\sR es}_\gO$, a natural structure of $\sO$-module. 
\bigskip

\begin{dfn}\label{nusector}{\bf\protect\boldmath $\nu$-sectorial region or $\nu$-sector\protect\unboldmath}\\
{\rm Given $\nu>0$, smaller than half the minimal distance between 
 the elements of $\gO$, we call  $\nu$-{\rm sectorial 
region} or $\nu$-{\rm sector} $\gD_\nu$ a domain of the Riemann surface $\sR_{\gO}$ composed of the three following parts:} 
{\rm
\begin{itemize}
\item an open sector $\gS_\nu$ with bounded opening at infinity;
\item a neighborhood of 0, say, an open disc $D_\nu$ centered at 0;
\item a tubular neighborhood $\sN_\nu$ of a  piecewise-$\sC^1$ path $\gamma$
connecting $D_\nu$ to $\gS_\nu$ after a finite number of turns
around  all or part of points of $\gO$.
\end{itemize}}

{\rm Moreover, the distance of $D_\nu$ to $\gO\setminus \{0\}$ and the distance of $\sN_\nu\cup \gS_\nu$ to $\gO$ has to be greater than $\nu$.}
\begin{center}
\includegraphics[scale=1]{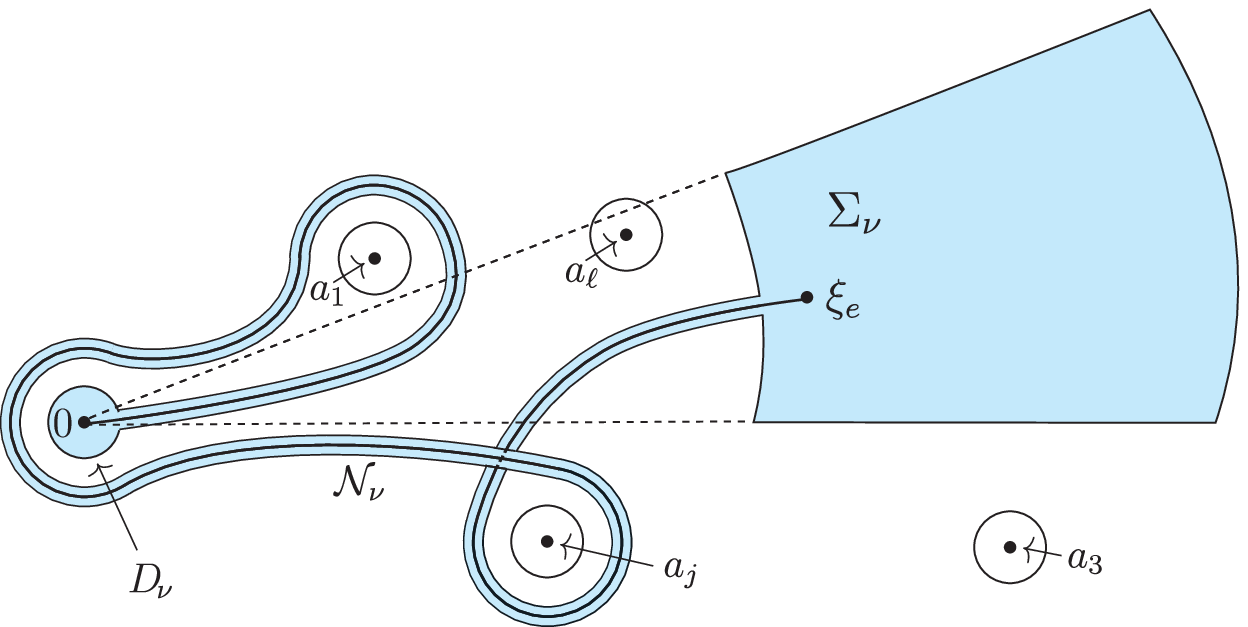}
\\ \quad
\end{center}
\end{dfn}

\begin{lemma}\label{nu-sector} Let $\gD_\nu$ be a sectorial region.

There exits a constant $K>0$ so that, for all $\xi\in \gD_\nu$, 
there is a piecewise-$\sC^1$-path $\gc_\xi$ contained in $\gD_\nu$ and parameterized 
by arc length from 0 to $\xi$ such  
that the arc length $s_\eta$ of all $\eta\in\gc_\xi$ satisfies
\begin{equation}\label{sxi}
\vert\eta\vert\leq s_\eta\leq K\vert \eta\vert.
\end{equation}
$\vert\eta\vert$ denotes the modulus of (the projection of) $\eta$ in $\C$.
\end{lemma}
\proof \ The first inequality being trivial we just have to prove that $s_\eta\leq K \vert \eta\vert$.

Assume first that the path $\gamma$ is $\sC^1$.\\
Let $\gamma'$ denote the extension of $\gamma$ from 0 to its beginning point, say, by a straight line. Let $\xi_e$ denote the end-point of $\gamma$ in $\gS_\nu$ and $\ell \vert \xi_e\vert$ denote the length of $\gamma'$ from 0 to $\xi_e$. The transversals of $\gamma$ in $\sN_\nu$ have all the same finite length.
\begin{itemize}
\item The property is clearly true in $D_\nu$ with $K=1$.
\item Given $\xi$ in $\sN_\nu$, consider a transversal issuing from $\xi$ and denote $\xi_0$ its intersection with $\gamma$. Let $\gamma_\xi$ be the path made of $\gamma'$ from 0 to $\xi_0$ followed by the arc of transversal $\displaystyle \mathop{\xi_0\xi}^\frown$.\\
For all $\eta\in \gamma_\xi$, the arc length $s_\eta$ is less than a constant $c$ independent of $\xi$ ($\gamma$ and the transversals have finite length)  and the euclidian distance $\vert\eta\vert$ is greater than $\nu$. Hence, $s_\eta\leq \frac{c}{\nu}\vert\eta\vert$.
\item Given $\xi\in\gS_\nu$, consider the point $\xi_0$ at distance $\vert\xi\vert$ from 0 on the ray issued from 0 through $\xi_e$. Let $\gamma_\xi$ be the path $\gamma'$ followed by the segment $(\xi_e\xi_0)$ and the arc $\displaystyle \mathop{\xi_0\xi}^\frown$ of the circle centered at 0 with radius $\vert \xi\vert$ included in $\gS_\nu$.\\
For $\eta\in(\xi_e\xi_0)$ the inequality holds with $K=\ell$. For $\eta $ belonging to the arc $\displaystyle \mathop{\xi_0\xi}^\frown$ we deduce from the inequalities 
$s_{\xi_0}\leq \ell\vert\xi_0\vert=\ell\vert\eta\vert$ and
$s_\eta-s_{\xi_0}\leq \alpha \vert \eta\vert$, where $\alpha$ denotes the opening angle of $\gS_\nu$, that the inequality holds with $K=\ell+\alpha$.
\end{itemize}
The case when the path is piecewise-$\sC^1$ is similar and left to the reader. \qed
\begin{dfn}\label{sum-res} {\bf Summable-resurgence.} 
{\rm
\begin{itemize}
\item A resurgent function $\hf(\xi)\in \widehat{\sR es}_{\gO}$ is
  said to be {\it summable-resurgent} when it grows at most
  exponentially at infinity on any $\nu$-sectorial region $\Delta_\nu$
  of $\sR_{\gO}$. 
\item A resurgent series $\tf(x)\in\widetilde{\sR es}_{\gO}$ is said to be {\it   summable-resurgent} if its Borel transform is a summable-resurgent  function. 
\end{itemize}
We denote respectively $\widehat{\sR es}_{\gO}^\text{\rm sum}$ and 
$\widetilde{\sR es}_{\gO}^\text{\rm sum}$ the set of summable-resurgent functions and the set of summable-resurgent series with singular support $\gO,0$. We denote $\widehat{\sR es}_{\gO,\go}^\text{\rm sum}$ the set of functions of $\widehat{\sR es}_{\gO}^\text{\rm sum}$ translated by $\go$.
}\end{dfn}

In general, the exponential type depends on the $\nu$-sectorial region
$\Delta_\nu$ and is unbounded when the width of $\Delta_\nu$ goes to
infinity (\cf Remark \ref{largeur}). 
\bigskip

The set  $\widehat{\sR es}_{\gO,\go}^\text{\rm sum}$ is a $\sO^{\leq
  1}(\C)*$submodule of $\widehat{\sR es}_{\gO,\go}$. 
\bigskip

\begin{rem}{\rm 
A summable-resurgent series is both resurgent
and summable but} {\it the converse is false}.
{\rm Indeed, a series is summable when it satisfies the conditions of Definition \ref{sum-res} in restriction to the first sheet only.}
\end{rem}

\subsection{Summable-resurgence theorem}

We are now able to state the result  in view in this section:
\begin{thm}  (J. \'Ecalle \cite{Ec85})\label{thmsum-res} \\
Assume that System (1) has a single level equal to one (\cf Assumption
(\ref{nivpur1})) and denote $\tf(x)$ the first $n_1$ columns of $\tF(x)$.

Then, $\tf$ is summable-resurgent with singular support ${\gO},0$ (recall $a_1=0$):
\renewcommand{\arraystretch}{2} 
\begin{center}
\begin{tabular}{| p{3cm}|}
\hline 
$\  \ \tf(x)\in \widetilde{\sR es}_{\gO}^\text{\rm sum}$\\
\hline
\end{tabular}\, .
\end{center}
\end{thm}

By means of elementary transformations (\cf Section \ref{prepared}) the result can be extended to the full matrix $\tF$ replacing however $\gO$ by ${\bf \gO}$. Therefore, we can state:
\begin{cor}
Under the same conditions as in Theorem \ref{thmsum-res} the full matrix $\tF$ is summable-resurgent with singular support ${\bf \gO},0$:
\renewcommand{\arraystretch}{2} 
\begin{center}
\begin{tabular}{| p{3cm}|}
\hline 
$\  \ \tF(x)\in \widetilde{\sR es}_{\bf\gO}^\text{\rm sum}$\\
\hline
\end{tabular}\, .
\end{center}
\end{cor}

The hypothesis of ``single level equal to 1'' is central. 
In a paper to come it will be shown how to generalize this result to the
case of a single level not equal to one. In the case of several levels
one has to consider several Borel planes simultaneously. 
\bigskip

Theorem \ref{thmsum-res} can be proved in different ways. 
We first  sketch a proof based on the Newton polygon and Ramis Index Theorem,  the system being given in the 
form of an equation of order $n$. Next, we develop a proof 
following \'Ecalle's method by regular perturbation of the system 
and majorant series. This second approach will allow us to 
precisely describe the singularities of the Borel transform $\gvf(\xi)=\hf(\xi)$ in the Borel plane which are all located in the set $\gO$. 


\subsection{Sketch proof of summable-resurgence using linear 
differential equations}\label{sumreseq}
The formal Borel transformation $\widetilde\sB$ is an isomorphism 
from the  differential algebra 
$\displaystyle \Big(\Cf, +,{\scriptscriptstyle\bullet},
x^2\frac{d}{dx}\Big)$ 
to the differential algebra 
$\displaystyle \Big(\gd\oplus\C\lbrack\lbrack\xi\rbrack\rbrack, +,
 \star, \xi{\scriptscriptstyle\bullet}\Big)$ that changes ordinary 
product into convolution product and changes derivation 
$\displaystyle x^2\frac{d}{dx}$ into multiplication  by $\xi$. 
It also  changes multiplication by $\displaystyle \frac{1}{x}$ 
into  derivation $\displaystyle \frac{d}{d\xi}$ allowing thus 
to extend the isomorphism from the meromorphic series 
$\Cf\lbrack1/x\rbrack$  to 
$\C\lbrack \gd^{(k)}, k\in \N\rbrack \oplus\C\lbrack\lbrack\xi\rbrack\rbrack$. Recall that the Borel transform of a monomial reads 
$\displaystyle \sB(x^m)=\frac{\xi^{m-1}}{\Gamma(m)}$ for all $m>0$.
\vspace{.5cm}

Consider now a differential equation $Dy(x)=0$ with  single level
one. By the Birkhoff Algebraization Theorem \cite[Th. 3.3.1]{Sib90} 
we may assume that the
operator $D$ has rational coefficients.  By means of an elementary, possibly trivial, gauge transformation  (\cf Section \ref{prepared}) we may also assume that the equation admits a formal series  solution $\tf(x)$. 

Multiplying $D$ by a convenient power of $1/x$ 
if needed, the Borel transformed equation $\widehat D \hy(\xi)=0$ 
 is again an ordinary linear differential equation 
with polynomial coefficients. Ramis Index Theorem \cite{R84} shows 
that the series  $\tf(x)$ is of Gevrey type of order 1. Hence, its 
Borel transform $\widetilde\sB(\tf)(\xi)$ converges in a neighborhood of the 
origin $\xi=0$;  we denote  by $\gvf(\xi)$ its sum.
A direct calculation using the characteristic equation associated 
with slope 1 of the Newton polygon of $D$ shows that the 
singularities of $\widehat D \hy(\xi)=0$ belong  to the 
finite set $\gO$ of the Stokes values of $Dy(x)=0$. From the 
Cauchy-Lipschitz theorem one can then assert that $\gvf(\xi)$ is 
resurgent with singular support $\gO,0$, 
\ie it can be analytically continued along any path issuing from 0 and staying in $\C\setminus \gO$. To see that  $\gvf(\xi)$  has exponential 
growth at infinity, it suffices to notice  that the equation 
$\tvi{14}{0}\widehat D \hy(\xi)=0$ has rank one at infinity 
\cite[Theorem 1.4]{Mal91}.

\subsection{Proof of summable-resurgence on systems following J.
  \'Ecalle's approach}\label{sumressystem} 
Since, by the Cyclic Vector Theorem, equations and systems 
are meromorphically equivalent, to get the summable-resurgence 
of solutions of systems there is no need for a new proof made 
directly on systems. However, the proof below quoted by J.~\'Ecalle 
\cite{Ec85} has its own interest. In particular, it allows us to give a precise description of the singularities of $\gvf(\xi)$ in the Borel plane.
A simpler case where the formal monodromy is assumed to be trivial ($L=I_n$) can be found in \cite{L-R95}; the case of a higher level can be found in \cite{Re07}.

\subsubsection{Setting the problem.}

 The conclusion of Theorem
\ref{thmsum-res} being preserved by elementary and meromorphic gauge transformations,
we can assume that the system is so prepared that  
conditions (\ref{F=I+}), (\ref{Re}) and (\ref{a1-0})  are all
satisfied.

The normal form (\ref{A0}) of System (\ref{A}) reads 
\begin{equation}\label{A00}
x^2\frac{dY}{dx}=A_0(x) Y \quad \text{ where }\quad
 A_0(x)=\bigoplus_{j=1}^J a_jI_{n_j}+x(\gl_jI_{n_j}+J_{n_j})
\end{equation}
while the matrix $A(x)$ of System  (\ref{A}) reads
\begin{equation}
A(x)=A_0(x)+B(x) \quad \text{ where }\quad B(0)=0
\end{equation}
More precisely, split  the matrix $B(x)=\begin{bmatrix}
B^{j;k}(x)
\end{bmatrix}$ into blocks fitting to the Jordan structure of
$L$. 
Then,
\begin{equation}\label{OB}
B^{j;k}(x)=\left\lbrace\begin{array}{lll}
O(x)& \text{if} & a_j\neq a_k\\
\noalign{\smallskip}
O(x^2)& \text{if} & a_j=a_k \text{ (and especially if } j=k).
\end{array}\right.
\end{equation}

As previously said, we restrict the study to the first 
column-block $\tf$ of $\tF$ uniquely determined by the first $n_1$ 
columns of the homological system
\begin{equation}\label{AA0}
x^2\frac{d\tF}{dx}-A_0(x)\tF+\tF A_0(x)=B(x)\tF
\end{equation}
jointly with the initial condition $\tf (0)=I_{n,n_1}$ (Recall that $I_{n,n_1}$ denotes  the identity matrix $I_n$ truncated at its first $n_1$ columns).
Hence, the system 
\begin{equation}\label{AA0f}
x^2\frac{d\tf}{dx}-A_0(x)\tf+x\tf J_{n_1}= B(x)\tf.
\end{equation}
(Recall $a_1=\gl_1=0$).

It is not clear from the Borel transformed system 
\begin{equation}\label{BB_0}
\xi \hf-\widehat A_0*\hf+1*\hf J_{n_1}=\widehat B*\hf
\end{equation}
that the Borel transform $\hf(\xi)$ of $\tf(x)$ satisfies the conditions of Definition \ref{sum-res} for $\tf(x)$ to be summable-resurgent.
\vspace{.5cm}

To prove the summable-resurgence of $\tf(x)$, J.~\'Ecalle suggests in
\cite{Ec85} to regularly perturb System (\ref{AA0f}) by substituting $\ga
B$ for $B$, next to solve this perturbed system  in terms of a power series
in the parameter $\ga$ and then, to proceed by majorant series satisfying a convenient system. There
exists, of course,
many possible majorant systems. Here below, we make explicit a 
possible one. We consider the very general case when $L$ is in non
diagonal form, covering thus, all possible cases of resonances. The calculation is made more complicated than
in the diagonal case since we have to work with packs of equations
instead of individual ones but the philosophy keeps the same.
\vspace{.5cm}

We split $\tf(x)$ into row-blocks $\tfjp{j}(x)$ accordingly to 
the Jordan structure of $L$ and we refer to Section \ref{knownresults}  for the notations.

\subsubsection{The perturbed system.}

An identification of equal powers of $\ga$ shows that the 
perturbed system
\begin{equation}\label{AA0falpha}
x^2\frac{d\tf}{dx}-A_0(x)\tf+x\tf J_{n_1}=\ga B(x)\tf
\end{equation} equivalent, for all $j\in \{1,\dots,J\}$, to
\begin{equation}\label{AA0jfalpha}
x^2\frac{d\tfjp{j}}{dx}-(a_j+\gl_j x)\tfjp{j} -xJ_{n_j}\tfjp{j}+x\tfjp{j}J_{n_1}=\alpha B^{j;\scriptscriptstyle{\bullet}}\tf
\end{equation}
 admits a unique formal solution of the form $\displaystyle 
\tf(x,\ga)=\sum_{m\geq 0}\tf_m(x)\ga^m$ satisfying
$\displaystyle 
\tf_0(x)=I_{n,n_1}$ and $\displaystyle  \tf_m(x)\in x\Cf $ 
for all $ m\geq 1$.
Since there is no possible ambiguity we keep denoting $\tf$ the 
perturbed solution. The proof proceeds by induction on $m$.  For all 
$j\in\{1,\dots,J\}$, System (\ref{AA0jfalpha})
reads, for $\ell=1,\dots, n_j$ and $k=1,\dots,n_1$,
\begin{equation}\label{AA0fm}
x^2\frac{d\tfjq{(j,\ell)}{(k)}_m}{dx}-(a_j+\gl_j x)\tfjq{(j,\ell)}{(k)}_m= 
x\tfjq{(j,\ell+1)}{(k)}_m -x\tfjq{(j,\ell)}{(k-1)}_m+ 
B^{(j,\ell);\scriptscriptstyle{\bullet}}\tfpq{(k)}_{m-1}
\end{equation}
and can be solved term after term following the alphabetic order on 
$(k,n_j-\ell)$ (one begins with the first column $k=1$ and  the last row $\ell=n_j$
in each block $j$).

It turns out actually that for all $j$ and $m\geq 1$,
$$\displaystyle \tfjp{j}_{2m-1}=O(x^{m}) \quad\text{ and } \quad \tfjp{j}_{2m}=
\left\lbrace\begin{array}{ll}
O(x^{m})& \text{ if } a_j=0\\
\noalign{\smallskip}
O(x^{m+1})& \text{ if } a_j\neq 0\\
\end{array}\right.  $$
allowing thus to rewrite the series
$\displaystyle 
\tf(x,\ga)=\sum_{m\geq 0}\tf_m(x)\ga^m$ in $\ga$ as a series
 in $x$ with coefficients that are polynomial in $\ga$. Consequently, the unperturbed solution $\tf(x)$ 
corresponds to $\tf(x,1)$ 
(unicity of  $\tf(x)$ and $\tf(x,1)$) and,  for all $\ga$ 
and in particular for $\ga=1$, the series $\displaystyle \tf(x,\ga)$ 
admits a formal Borel transform $\hf(\xi,\ga)$ with respect to $x$ 
of the form
\begin{equation}
\hf(\xi,\ga)=\gd I_{n,n_1}+\sum_{m\geq 1}\gvf_m(\xi)\ga^m
\end{equation} 
where $\gvf_m(\xi)$ denotes the Borel transform of $\tf_m(x)$ for 
all $m\geq1$. 

For fixed $j\in \{1,\dots,J\}$ and $m\geq 1$ the system satisfied 
by the $\gvf^{(j,\ell);(k)}_m(\xi)$'s reads, for $\ell=1,\dots,n_j$ 
and $k=1,\dots,n_1$,
\begin{equation}\label{systemphim}
\left\lbrace\begin{array}{l}\displaystyle 
(\xi-a_j)\frac{d\gvf_m^{(j,\ell);(k)}}{d\xi}-(\gl_j-1)\gvf_m^{(j,\ell);(k)}=\\
\noalign{\medskip}
\displaystyle \hspace{1.5cm} \gvf_m^{(j,\ell+1);q}-\gvf_m^{(j,\ell);(k-1)}+ 
\frac{d\hBjp{(j,\ell)}}{d\xi}*\gvf_{m-1}^{\scriptscriptstyle\bullet;(k)}
+\hBjp{(j,\ell)}(0) \gvf_{m-1}^{\scriptscriptstyle\bullet;(k)} 
\end{array}\right.
\end{equation}
since the Borel transform of $\displaystyle \frac{1}{x} B^{(j,\ell);\scriptscriptstyle{\bullet}}(x)\tfpq{(k)}_{m-1}(x)$ is equal to 
$$
\frac{d}{d\xi}\Big( \hBjp{(j,\ell)}*\gvf_{m-1}^{\scriptscriptstyle\bullet;(k)} \Big)(\xi) =  \frac{d\hBjp{(j,\ell)}}{d\xi}*\gvf_{m-1}^{\scriptscriptstyle\bullet;(k)} +
\hBjp{(j,\ell)}(0) \gvf_{m-1}^{\scriptscriptstyle\bullet;(k)}
$$
\noindent ($\hBjp{(j,\ell)}(0)$ is also the coefficient of $x$ in $B^{(j,\ell);\scriptscriptstyle\bullet}(x)$ and is not supposed equal to zero).\\
It can again be solved term after term starting with the 
first column $k=1$ and  the last row  $\ell=n_j$.
\medskip

Recall that for $b$ an entire function, 
the convolution $b*\xi^m, m\in\N$ is well defined by the integral 
$\int_0^\xi b(\xi-t)t^m dt$ as well as $b*\hf$ by 
$\int_\gc b(\xi-t)\hf(t) dt$ along any path $\gc$ avoiding the singularities of $\hf$ (here, the Borel transforms  
$\displaystyle\hBjp{(j,k)} $ of $B^{(j,k);\scriptscriptstyle{\bullet}}$  
are entire functions since $B(x)$ is analytic at 0).
\medskip

Note also that the only  singularities 
of System (\ref{systemphim}) are the Stokes values $a_j\in \gO$ and we can conclude that the $\gvf_m$'s are resurgent functions defined on $\sR_\gO$.

\subsubsection{What has to be proved. }

Fix now $\nu>0$ and a $\nu$-sectorial region 
$\gD_\nu$ as described in Definition \ref{nusector}. 

We are left to prove that
\begin{description}
\item{(a)} the Borel series 
 $\hf(\xi)=\hf(\xi,1)$ is convergent and can be analytically 
continued to $\sR_{\gO}$; we keep denoting $\hf(\xi)$ the analytic continuation;
\item{(b)} $\hf(\xi)$ grows at most exponentially on $\gD_\nu$ at infinity.
\end{description}

These properties could be directly shown to be true for the $\gvf_m$'s. To prove that they are true for $\hf$ we  use the technique of majorant series.

\subsubsection{A candidate majorant system.}

Instead of System (\ref{AA0jfalpha}) consider, for $j=1,\dots,J$, 
the perturbed linear system
\begin{equation}\label{systemg}
\left\lbrace 
\begin{array}{l}
\displaystyle 
C_j(\tgjp{j}- \hspace{-2pt}I_{n,n_1}^{j;\scriptscriptstyle{\bullet}})=\displaystyle 
J_{n_j}\tgjp{j}+\tgjp{j}J_{n_1} \hspace{-2pt}
-2I_{n,n_1}^{j;\scriptscriptstyle{\bullet}} J_{n_1}+\ga \frac{\vert
  B\vert^{j;\scriptscriptstyle{\bullet}}}{x}\tg  \text{  if } a_j= 0\\
 \noalign{\bigskip}
(\nu-\vert \gl_j-1\vert x)\tgjp{j}=J_{n_j}x\tgjp{j}+x\tgjp{j}J_{n_1} 
 +\ga\vert B\vert^{j;\scriptscriptstyle{\bullet}} \tg\qquad   \text{  if } a_j\neq 0 \ \\
\end{array}
\right.
\end{equation}
where the unknown $\tg$ is, like $\tf$, a $n\times n_1$-matrix split into row-blocks $\tgjp{j}$ fitting the Jordan structure of $L$  
and where $\vert B\vert$ denotes the series $B$ in
which the coefficients of the powers of $x$ are replaced by their module. The constants $C_j>0$ are to be  adequately chosen which we will do in Lemma \ref{lemmeseriemaj} below. For now, they are just arbitrary non-zero constants.

\medskip

System (\ref{systemg}), like System  (\ref{AA0jfalpha}), admits a unique formal 
solution $\displaystyle \tg(x,\ga)=\sum_{m\geq 0}
\tg_m(x)\ga^m$ such that $\tg_0(x)=I_{n,n_1}$ and $\tg_m(x)\in x\Cf$ for all
$m\geq 1$. Like $\tf(x,\ga)$, this series  satisfies 
$$\displaystyle \tgjp{j}_{2m-1}=O(x^{m}) \quad\text{ and } \quad \tgjp{j}_{2m}=
\left\lbrace\begin{array}{ll}
O(x^{m})& \text{ if } a_j=0\\
\noalign{\smallskip}
O(x^{m+1})& \text{ if } a_j\neq 0\\
\end{array}\right. $$
and consequently, $\tg(x,\alpha)$ can be seen as a series in powers of $x$ whose coefficients are polynomials in $\alpha$.
But unlike $\tf_m(x)$ in general, $\tg_m(x)$ has 
non negative coefficients for all $m\geq 1$.
\medskip

Prove now that $\tg(x,\ga)$ is a convergent series in $(x,\ga)$ in a domain containing $(x,\ga)=(0,1)$.

From System (\ref{systemg})  written for each individual column $k$ we obtain, for all $j$, the linear system 
\begin{equation}\label{systemgq}
\left\lbrace
\begin{array}{ll}\displaystyle 
(C_j-J_{n_j})\tg^{j;(k)} - \alpha \frac{\vert B\vert^{j;\scriptscriptstyle\bullet}}{x} \tg^{{\scriptscriptstyle\bullet};(k)} =\tg^{j;(k-1)}+ {\rm const}  \quad&\text{ if }a_j=0\\
\noalign{\medskip}
\Big(\nu-\vert\gl_j-1\vert x-J_{n_j}x\Big)\tg^{j;(k)}-\alpha\vert B\vert^{j;\scriptscriptstyle\bullet} \tg^{{\scriptscriptstyle\bullet};(k)} =x\tg^{j;(k-1)} & \text{ if } a_j\neq 0.\\
\end{array}\right.
\end{equation}
\noindent (We set $\tg^{j;(0)}=0$).\\
For $x=0$ the system reduces to 
\begin{equation}\label{systemgqx=0}
\left\lbrace
\begin{array}{ll}\displaystyle 
(C_j-J_{n_j})\tg^{j;(k)} - \alpha \vert B'\vert^{j;\scriptscriptstyle\bullet}(0) \tg^{{\scriptscriptstyle\bullet};(k)} =\tg^{j;(k-1)}+ {\rm const}  \quad&\text{ if }a_j=0\\
\noalign{\medskip}
\nu \tg^{j;(k)} =0 & \text{ if } a_j\neq 0.
\end{array}\right.
\end{equation}
However, $\displaystyle \vert B'\vert^{j;\scriptscriptstyle\bullet}(0) \tg^{{\scriptscriptstyle\bullet};(k)}=\sum_{r\,|\,a_r\neq 0} \vert B'\vert^{j;r}(0)\tg^{r;(k)}$ since $\vert B'\vert^{j;r}(0)=0$ for all $r$ such that $a_r=a_j=0$ and an adequate linear combination among equations of this system allows to cancel the terms in $\alpha$. 
System (\ref{systemgqx=0}) is thus equivalent to a constant triangular system whose diagonal terms are either equal to $C_j$ or to $\nu$. Having assumed $C_j\neq 0$ for all $j$ it is then a Cramer system and System (\ref{systemgq}) is equivalent to a system of the form 
$\displaystyle 
(T+xM(x,\ga))g^{\scriptscriptstyle{\bullet};(k)}=N(g^{\scriptscriptstyle{\bullet};(k-1)})$
where $T$ is a constant  invertible matrix and $M$ is analytic in $(x,\alpha)$ in a strip around $x=0$  while the right hand-side 
is a column vector depending on  $g^{j;(k-1)}$ and analytic in the same strip. 
The determinant of System (\ref{systemgq}) is analytic in $(x,\alpha)$ and non-zero for $x=0$ whatever $\alpha$ is equal to. There exists then a bi-disc centered at $(x,\alpha)=(0,0)$ and containing $(x,\alpha)=(0,1)$ on which the determinant does not vanish. On such a bi-disc, 
System (\ref{systemgq}) admits an analytic solution $g(x,\alpha)$. 
By unicity, its Taylor expansion $g(x,\alpha)=\sum g_m(x) \alpha^m$ coincides with $\tg(x,\alpha)=\sum \tg_m(x) \alpha^m$. 

In particular, for $\alpha=1$, the formal solution $\tg(x)=\tg(x,1)=\sum\tg_m(x)$ converges.  Henceforth,  its Borel transform $\hg(\xi)=\sum\phi_m(\xi)$ is an entire function with exponential growth at infinity ($\phi_m$ denotes the Borel transform of the series $\tg_m(x)$ and is also an entire function with exponential growth at infinity).

\subsubsection{Majorant series and exponential growth.}

Lemma \ref{lemmeseriemaj}  below shows that $\hg(K\vert\xi\vert)=\sum \phi_m(K\vert\xi\vert)$ is a majorant series 
for $\hf(\xi)=\sum \gvf_m(\xi)$ on $\Delta_\nu$.
\vspace{.5cm}

\begin{lemma}\label{lemmeseriemaj} \quad

Let $K>0$ be associated with the chosen $\nu$-sector $\Delta_\nu$ as in Lemma \ref{nu-sector}.

For all $m\geq 0,\ \xi
  \in\gD_\nu,\ j=1,\dots,J$ and $q=1,\dots,n_1$, the following inequalities hold:
\begin{equation}
\vert \gvf_m^{j;(k)}(\xi)\vert \leq \phi_m^{j;(k)}(s_\xi)
 \leq \phi_m^{j;(k)}(K\vert\xi\vert).
\end{equation}
($\gvf_m$ and $\phi_m$ are the Borel transforms of the ``initial'' solution $\tf_m$ and of the ``majorant'' solution $\tg_m$ respectively).
\end{lemma}
\proof Recall that the functions $\gvf_m$ are defined over $\sR_\gO$ and the functions $\phi_m$ on all of $\C$. They are then all well defined over $\Delta_\nu$.

$\bullet$ The series $\phi_m(\xi)$ have non negative coefficients and by Lemma \ref{nu-sector} we know that we can connect $\xi$ to 0 by a path so that $s_\xi\leq K\vert\xi\vert$. Hence, the second inequality.

$\bullet$ Prove the first inequality.
For all $m\geq 1$, the entries $\gvf_m^{(j,\ell);(k)}$ (row $\ell$ of 
row-block $j$ and column
$k$) of $\gvf_m$ and $\phi^{(j,\ell);(k)}_m$ of $\phi_m$  satisfy respectively 
\begin{equation}\label{gvfm}
(\xi-a_j)\frac{d}{d\xi}\gvf_m^{(j,\ell);(k)}-(\gl_j-1)\gvf_m^{(j,\ell);(k)} =\go^{(j,\ell);(k)}_m
\end{equation}
where
\begin{equation*}
\go^{(j,\ell);(k)}_m=\gvf^{(j,\ell+1);(k)}_m-\gvf_m^{(j,\ell);(k+1)}+\frac{d}{d\xi}
\hBjp{(j,\ell)} *\gvfpq{(k)}_{m-1} 
+\hBjp{(j,\ell)}(0) \gvf_{m-1}^{\scriptscriptstyle\bullet;(k)}
\end{equation*}
 and 
\begin{equation}
\left\lbrace\begin{array}{cccl}\displaystyle 
\nu\frac{d\phi^{(j,\ell);(k)}_m}{d\xi}-\vert \gl_j-1\vert\phi^{(j,\ell);(k)}_m
&=&\gO^{(j,\ell);(k)}_m
&\quad\text{if } a_j\neq 0\\
 \noalign{\bigskip}\displaystyle 
 C_j\phi^{(j,\ell);(k)}_m&=&\gO^{(j,\ell);(k)}_m
&\quad\text{if } a_j= 0\\
\end{array}\right.
\end{equation}
where 
\begin{equation*}\displaystyle 
\gO^{(j,\ell);(k)}_m=\phi^{(j,\ell+1);(k)}_m+\phi^{(j,\ell);(k-1)}_m
+ \frac{d}{d\xi}\vBjp{(j,\ell)}*\phipq{(k)}_{m-1}
+\vBjp{(j,\ell)}(0) \phi_{m-1}^{\scriptscriptstyle\bullet;(k)}.
\end{equation*}

Fix $\xi\in\gD_\nu$ and a path $\gc_{\xi}$ in $\gD_\nu$ as in Lemma \ref{nu-sector} so that  $$\vert \eta\vert \leq s_\eta
\leq K\vert\eta\vert \quad\text{ for all } \eta\in\gc_{\xi}.$$
We proceed by recurrence following the alphabetic order 
on $(m,k,n_j-\ell)$ and we assume that for all 
$(m',k',n_j-\ell')<(m,k,n_j-\ell)$, the inequality
$\vert \gvf_{m'}^{(j,\ell');(k')}(\eta)\vert \leq
\phi_{m'}^{(j,\ell');(k')}(s_\eta)$
holds for all $\eta\in\gc_{\xi}$. (It holds for $m'=0$).\\
We observe that 
$\vert \go^{(j,\ell);(k)}_m(\xi)\vert\leq \gO^{(j,\ell);(k)}_m(s_\xi)$.
Indeed,
$$
\Big\vert \frac{d\hBjp{(j,\ell)}}{d\xi} *\gvfpq{(k)}_{m-1}(\xi)\Big \vert=
\displaystyle\Big \vert \int_0^{s_\xi}
\frac{d\hBjp{(j,\ell)}}{d\xi}(\xi-\eta(s))
\gvfpq{(k)}_{m-1}(\eta(s))\eta'(s) ds\Big\vert\qquad
$$
$$
\begin{array}{ccl}\displaystyle 
\qquad&\leq&\displaystyle \int_0^{s_\xi}
\frac{d\vBjp{(j,\ell)}}{d\xi}(\vert\xi-\eta(s)\vert)
\vert \gvfpq{(k)}_{m-1}(\eta(s))\vert ds\quad(\text{since }\vert\eta'(s)\vert=1)\\
\noalign{\bigskip}
&\leq& \displaystyle \int_0^{s_\xi}
\frac{d\vBjp{(j,\ell)}}{d\xi}(s_\xi-s)
\phipq{(k)}_{m-1}(s) ds\  \\
\noalign{\smallskip}
&&\qquad\qquad\qquad
(\text{since }\vBjp{(j,\ell)} \text{ has non negative coefficients})\\
\noalign{\smallskip}
&=&\displaystyle \frac{d\vBjp{(j,\ell)}}{d\xi} *\phipq{(k)}_{m-1}(s_\xi)\\
\end{array}
$$
and the other three terms of $\go^{(j,\ell);(k)}_m$ are majored using the recurrence hypothesis.
\medskip

 To conclude in the case when $a_j=0$ 
we solve Equation (\ref{gvfm}). Hence,
$$
\gvf^{(j,\ell);(k)}_m(\xi)=\xi^{\gl_j-1}\int_0^{s_\xi}
\go_m^{(j,\ell);(k)}(\eta(s))\eta(s)^{-\gl_j}\eta'(s)ds
$$
and then, $$ \vert \gvf^{(j,\ell);(k)}_m(\xi)\vert\hspace{-.2em}\leq\displaystyle 
\vert \xi\vert^{{\rm Re}\gl_j-1} e^{-{\rm Im}(\gl_j-1){\rm arg}(\xi)}\int_0^{s_\xi} \hspace{-.5em} \gO_m^{(j,\ell);(k)}(s)
\vert \eta(s)\vert^{-{\rm Re}\gl_j}e^{{\rm Im}(\gl_j){\rm arg}(\eta(s))}ds.$$
Since points in $\Delta_\nu$ have bounded arguments there exists a constant $c>0$ such that $\displaystyle \tvi{15}{0} e^{-{\rm Im}(\gl_j-1){\rm arg}(\xi)+ {\rm Im}(\gl_j){\rm arg}(\eta(s))}\leq c$ for all $\xi$ and $\eta(s) \in \Delta_\nu$. 
It results that 
\begin{equation*}
\begin{array}{ccl}
\vert \gvf^{(j,\ell);(k)}_m(\xi)\vert&\leq&\displaystyle 
c\vert \xi\vert^{{\rm Re}\gl_j-1} \int_0^{s_\xi}\gO_m^{(j,\ell);(k)}(s)
\vert \eta(s)\vert^{-{\rm Re}\gl_j}ds\\
\noalign{\bigskip}
&\leq& \displaystyle c\vert \xi\vert^{{\rm Re}\gl_j-1}
\gO_m^{(j,\ell);(k)}(s_\xi)\int_0^{s_\xi}\Big(\frac{s}{K}\Big)^{-{\rm Re}\gl_j}ds   \\
\noalign{\smallskip}
&&\qquad \text{ (by  Lemma  } \ref{nu-sector}, \text{ the fact that }  {\rm Re}\gl_j\geq 0     \\
&& \qquad\qquad
  \text{ and that }  \gO_m \text{ has   non  negative coefficients}) \\
\noalign{\bigskip}
&\leq& \displaystyle c\vert \xi\vert^{{\rm Re}\gl_j-1}
\frac{K\gO_m^{(j,\ell);(k)}(s_\xi)}{1-{\rm Re}\gl_j}\Big(\frac{s_\xi}{K}\Big)^{1-{\rm Re}\gl_j}\\
\end{array}
\end{equation*}
\begin{equation*}
\begin{array}{ccl}
\qquad\ &\leq& \displaystyle \frac{cK}{1-{\rm Re}\gl_j}\gO_m^{(j,\ell);(k)}(s_\xi)\\
&& \hfill \quad \text{ (using Lemma } \ref{nu-sector} 
\text{ and } 1-{\rm Re}\gl_j>0)\\
\noalign{\bigskip}
&=&\phi^{(j,\ell);(k)}_m(s_\xi) \text{ if we choose } \displaystyle C_j=\frac{1-{\rm Re}\gl_j}{cK}.\end{array}
\end{equation*}

To conclude in the case when $a_j\neq 0$ we apply Gr\"onwall Lemma to 
$\gvf^{(j,\ell);(k)}_m(\gc_\xi(s))$.
This achieves the proof of Lemma \ref{lemmeseriemaj}. \qed

Lemma \ref{lemmeseriemaj} shows that $\hg(K\vert\xi\vert)=\sum \phi_m(K\vert\xi\vert)$ is a majorant series 
for $\hf(\xi)=\sum \gvf_m(\xi)$ on $\Delta_\nu$.  Since $\hg$ is well defined on $\Delta_\nu$ with exponential growth at infinity the same property holds for $\hf(\xi)$ which achieves the proof of Theorem \ref{thmsum-res}. \qed

\begin{rem}\label{unifconv}
{\rm It results from the above proof that the series $\sum \gvf_m(\xi)$ converges uniformly to $\hf(\xi)$ on compact sets of $\sR_\gO$.}
\end{rem}

\begin{rem}\label{largeur}{\rm 
We see from Lemma \ref{lemmeseriemaj} that when $\hg$ grows 
exponentially with type
$a$, \ie satisfies an inequality $\vert\hg(\xi)\vert\leq
const. e^{a\vert\xi\vert}$ for large $\xi$ then $\hf(\xi)$ grows exponentially
with type $Ka$. When the width of the domain $\gD_\nu$ goes to infinity so does $K$; hence, the necessity for considering $\nu$-sectorial regions with bounded width. Also the estimates in the proof of Lemma \ref{lemmeseriemaj}  would no more be valid on sectors with unbounded width.}
\end{rem}
\begin{rem}{\rm 
One should think at reading the  previous two sections that the proof
with systems is much longer. This is not the case. If the same level
of detail were provided the proof on equations like sketched in
Section \ref{sumreseq} would be much longer and also less elementary
since it includes the Index Theorem for rank one and the Main
Asymptotic Existence Theorem at infinity.}
\end{rem}


\section{Singularities in the Borel plane}
 Theorem \ref{thmsum-res} tells us  that $\tf(x)$ is a resurgent series of $\widetilde{\sR es}_{\gO}^\text{\rm sum}$. Its Borel transform $\hf(\xi)$  is then, in particular, analytic on the Riemann surface $\sR_\gO$,  its possible singularities being the points $a_1, a_2,\dots,a_J$ of $\gO$ including $a_1=0$ out of the first sheet. 

The form of System (\ref{systemphim}) shows that the singularities of $\hf$ should at least involve poles since some $\gl_j$ are equal to 0, complex powers when some $\gl_j$ are not 0 and logarithms. It is already known that, in the case of a system with the unique level one,  the singularities all belong to the Nilsson class (\cf \cite{BJL81} for instance). In that case, the exponentials $\displaystyle e^{-a_j/x}$ in the formal fundamental solution $\tY(x)$ act as translations in the Borel plane; hence, the location of the singularities at the various points $a_j$. In the case of higher level or in the case of several levels that, after rank reduction, we could assume to be all $\leq 1$ there might also occur exponentials of degree less than one. Their action in the Borel plane would be then transcendental and would generate  irregular ---no longer in the Nilsson class---  singularities.
\medskip

 Our aim in this section is to set up a precise description of the  singularities  in the Borel plane  related to the form of the formal fundamental solution $\tY(x)=\tF(x) x^L e^{Q(1/x)}$ of System (\ref{A}). 
 \medskip

 As previously, we restrict our study to the first $n_1$ columns $\tf$ of $\tF$. We base the analysis on the results of Section 
 \ref{sumressystem} interpreting the Borel transform $\hf(\xi)$ of $\tf(x)$ as a series $\hf(\xi)=\gd I_{n,n_1} + \sum_{m\geq 1}\gvf_m(\xi)$ which converges uniformly on  compact sets of $\sR_\gO$.
 
 Decomposing  $\gvf_m$ into blocks $ \displaystyle \left\lbrack\gvf_m^{j,\scriptscriptstyle\bullet}\right\rbrack_{1\leq j\leq J}$
 fitting the Jordan structure of $L$ ($\gvf_m^{j,\scriptscriptstyle\bullet}$ has dimension $n_j\times n_1$), System (\ref{systemphim}) splits into the following $J$ systems:
\begin{equation}\label{systemphimblocs}
(\xi-a_j)\frac{d\gvf_m^{j,\scriptscriptstyle\bullet}}{d\xi}-(\gl_j-1)\gvf_m^{j,\scriptscriptstyle\bullet}-J_{n_j} \gvf_m^{j,\scriptscriptstyle\bullet}+\gvf_m^{j,\scriptscriptstyle\bullet}J_{n_1} =\frac{d\hBjp{j}}{d\xi}*\gvf_{m-1}
+\hBjp{j}(0) \gvf_{m-1}
\end{equation}
where  $\displaystyle \frac{d\hBjp{j}}{d\xi}*\gvf_{m-1} = \sum_{k=1}^J \frac{d\hB^{j;k}}{d\xi}*\gvf_{m-1}^{k;\scriptscriptstyle{\bullet}}$ and $\hBjp{j}(0) \gvf_{0}=0$.
\bigskip

Let us first introduce some vocabulary used in resurgence theory.
Working locally we place ourselves at the origin of $\C$. As previously, we denote $\sO=\C\{x\}$ the space of holomorphic germs at 0 in $\C$ and $\stO$ the space of holomorphic germs at 0 on the Riemann surface $\widetilde\C$ of the logarithm.

Being mostly interested in integrating solutions in the Borel plane on both sides of the singularities, thus enclosing them in a loop, we can neglect holomorphic terms and it is natural to consider the quotient space $\sC=\stO/\sO$. The elements of $\sC$ are called {\it micro-functions} by B. Malgrange \cite{Mal91} by analogy with hyper-  and micro-functions defined by Sato, Kawai and Kashiwara in higher dimensions. They are called {\it singularities} by J. Ecalle and al. and usually denoted with a nabla, like ${\buildrel {\scriptscriptstyle\,\nabla}\over \gvf}$,  for a singularity of the function $\gvf$ while the space $\sC=\stO/\sO$ is denoted ${\textsc{sing}}_0$ (\cf \cite{Sau05}). 
A representative of ${\buildrel {\scriptscriptstyle\,\nabla}\over \gvf}$ in $\stO$ is often denoted $\widecheck{\gvf}$ and is called {\it a major} of ${\gvf}$.\\
It is worth to consider the two natural maps 
\begin{equation*}
\begin{array}{rcccl}
{\rm can }:\quad &\stO& \longrightarrow & \sC=\stO/\sO& \text{ the canonical quotient map}\\
\noalign{\bigskip}
\text{and \ var }:\quad&\sC& \longrightarrow &\stO& \text{ the variation map,}\\
\end{array}
\end{equation*}
 action of a positive turn around 0 defined by $\displaystyle {\rm var }\, 
{\buildrel {\scriptscriptstyle\,\nabla}\over \gvf} (\xi)=
\widecheck{\gvf}(\xi)-\widecheck{\gvf}(\xi e^{-2\pi i})$ where $\widecheck{\gvf}(\xi e^{-2\pi i})$ is the analytic continuation of $\widecheck{\gvf}(\xi)$
along  a path turning once clockwise around 0 close enough to 0 for $\widecheck{\gvf}$ to be defined all along (the result is independent of the choice of a major $\widecheck{\gvf}$). The germ $\widehat\gvf={\rm var }{\buildrel {\scriptscriptstyle\,\nabla}\over \gvf}$ is called {\it the minor} of ${\buildrel {\scriptscriptstyle\,\nabla}\over \gvf}$.\\
Let $\delta$ denote the Dirac distribution at 0, $\delta^{(m)}$ its $m^{th}$ derivative and $Y=\partial_\xi^{-1}\delta$ the Heaviside (micro-)function. One can make the following identifications:
\begin{equation*}\begin{array}{ll}\displaystyle 
{\rm can}\Big(\frac{1}{2\pi i\xi}\Big)=\delta &\quad \displaystyle 
{\rm can}\Big(\frac{(-1)^m m!}{2\pi i\xi^{m+1}}\Big)=\delta^{(m)} \\
\noalign{\bigskip}
\displaystyle {\rm can}\Big(\frac{\ln\xi}{2\pi i}\Big)=Y \quad&\displaystyle 
{\rm can}\bigg(\Big(\frac{\ln\xi}{2\pi i}\Big)^2\bigg)=\Big(2\frac{\ln\xi}{2\pi i}-1\Big)Y\ \text{ and so on}\dots\\
\end{array}
\end{equation*}
It is sometimes useful not to work at the origin. 
Given $\go\neq 0$ in $\C$ we denote $\sC_\go=\textsc{sing}_\go$ the space of the singularities at $\go$, \ie the space $\sC=\textsc{sing}_0$ translated from 0 to $\go$.
A function $\widecheck\gvf$ is a major of  a singularity at $\go$ if $\widecheck\gvf(\go+\xi)$ is a major of a singularity at 0.

\subsection{Simple-moderate singularities}

In this Section, we state some properties, used further on, of the singularities ---poles, logarithms and complex powers--- 
which should occur in the Borel plane. We shall see (\cf Thm \ref{thmsingf}) that poles, logarithms and complex powers are  the only possible singularities arising in the Borel plane for linear systems with the unique level one, far from being any singularity in the Nilsson class.

\begin{dfn} {\bf Simple-moderate singularities}\label{footnotesimple}

{\rm \begin{itemize}
\item 
A singularity or micro-function ${\buildrel {\scriptscriptstyle\,\nabla}\over \gvf}$ at 0 is said to be {\it simple} if it has a major of the form\footnote{This definition of a simple singularity is less restrictive than the one one can find in the literature (\cf \cite{Ec85} or \cite{Sau05} for instance) where $N_0$ is taken equal to 0. Here, we allow powers of logarithms; still poles are required to be simple but they can be factored by logarithms. We will see that, in the linear case when the system is prepared like in Section \ref{prepared}, simple singularities in the restrictive sense would occur only under strong assumptions such as trivial  formal monodromy. For a general level-one system, not in prepared form, there could also occur multiple poles.}
\begin{equation*}
\widecheck{\gvf}(\xi)=\sum_{p=0}^{N_0} \ga_p \frac{\ln^p(\xi)}{\xi} +\sum_{p=1}^{N_0\,+1} \hh_p(\xi)\,\ln^p(\xi)
\end{equation*}
where $N_0\in\N$, $\ga_p\in\C$ and $\hh_p\in \sO$ for all $p$.
\item A singularity  or micro-function ${\buildrel {\scriptscriptstyle\,\nabla}\over \gvf}$ at 0 is said to be {\it simple-moderate} if it has a major which differs from a simple one by terms of  the form 
\begin{equation*}
\sum_{\gl\in\Lambda}\sum_{p=0}^{N_\gl} \ga_{\gl,p}\, \xi^{\gl-1}\,\ln^p(\xi) + \sum_{\gl\in\Lambda}\sum_{p=0}^{N_\gl}
\hH_{\gl,p}(\xi) \, \xi^{\gl}\, \ln^p(\xi)
\end{equation*}
where $\Lambda$ is a finite set of numbers $\gl\in\C$ satisfying $0<{\rm Re} \gl<1$ and for all $\gl\in\Lambda$, $N_\gl\in\N$, $\ga_{\gl,p}\in\C$ and $\hH_{\gl,p}\in \sO$ for all $p$.
\item A singularity or micro-function ${\buildrel {\scriptscriptstyle\,\nabla}\over \gvf}$ at $\go$ is said to be {\it simple} or {\it simple-moderate} if it has a major $\widecheck{\gvf}$ such that $\widecheck{\gvf}(\go+\xi)$ be of the previous forms respectively.
\end{itemize}
}\end{dfn}

Let $\gc$ be a path from 0 to $\xi$ in $\C$.\\
We denote $u*_\gc v(\xi)$ the convolution  product along the path $\gc$ defined by 
\begin{equation*}
u*_\gc v(\xi)=\int_\gc u(\xi-t)\, v(t)\, dt,
\end{equation*}
when the integral makes sense.

\begin{lemma}\label{convlog}{\bf Convolution with powers and logarithms}\\
Let $\gvf$ be an entire function on $\C$.\\
Let $p\in\N$ and $\lambda\in\C$ satisfying $0\leq {\rm Re}\,\lambda <1$.\\
Let $\go\in \C^*$, let $\psi$ be a function satisfying $\psi(\go)\neq 0$ and holomorphic on a domain containing 0 and $\go$ and let $\gc_\xi$ be a path from 0 to  $\xi$ avoiding $\go$ and contained in this domain.\\

Then, the convolution product $\gvf*_{\gc_\xi}\Big((\xi-\go)^{\lambda-1} \ln^p(\xi-\go)\psi\Big)(\xi)$ exists and is, close to $\go$, of the form
\begin{equation*}
\gvf*_{\gc_\xi}\Big((\xi-\go)^{\lambda-1} \ln^p(\xi-\go)\psi\Big)(\xi)=(\xi-\go)^\gl P\big((\ln(\xi-\go)\big) +{\rm ent} (\xi)
\end{equation*}
where $P(X)\in\C\{\xi\}[X]$ is a polynomial with holomorphic coefficients at $\go$ and degree 
\begin{equation*}
{\rm deg} (P)=\left\lbrace\begin{array}{lll}
p& {\rm if}& \gl\neq 0\\
p+1& {\rm if}& \gl= 0\\
\end{array}\right.
\end{equation*}
and where {\rm ent} stands for an entire function.
\end{lemma}
Note that the power of $(\xi-\go)$ in the right-hand side has increased by 1 unit at least.

\proof The convolution product  $\gvf*_{\gc_\xi}\Big((\xi-\go)^{\lambda-1} \ln^p(\xi-\go)\psi\Big)(\xi)$ is well-defined by the integral
$\displaystyle \int_{\gc_\xi} \gvf(\xi-t) (t-\go)^{\gl-1} \ln^p(t-\go)
\psi(t)dt$ since both $\gvf(\xi-t)$ and $\psi(t)$ are holomorphic along $\gc_\xi$; and in the case when $\go=0$, the factor $t^{\gl-1}$ is integrable at 0.
\medskip

We are interested in the behavior  of this function as $\xi$ goes to $\go$. \medskip


\noindent\begin{minipage}[c]{80mm}
\baselineskip=1.5em

Suppose $\xi$ so close to $\go$ that there is a disc $D_\go$ centered at $\go$, containing $\xi$ and included in the holomorphy domain of $\psi$. \\
Fix $a$ on $\gc_\xi$ so that  the part ${}_a\hspace{-.2em}\gc_\xi$ of $\gc_\xi$ from $a$ to $\xi$ belong to $D_\go$ and be homotopic to a straight line in $D_\go\setminus \lbrace\go\rbrace$.\\

\end{minipage}
\hskip .5cm
\begin{minipage}[c]{50mm}
\includegraphics[scale=1]{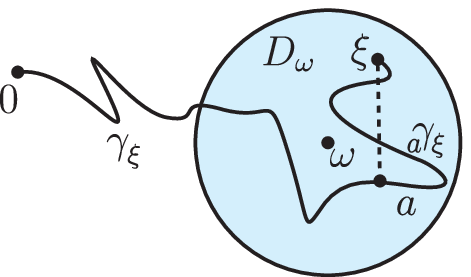}
\end{minipage}

Up to an entire function, we can replace the convolution product under consideration by the integral
 \begin{equation*}
g_a(\xi)=\int_a^\xi \gvf(\xi-t)\,(t-\go)^{\gl-1}\,\ln^p(t-\go)\,\psi(t)dt
\end{equation*}
in $\C$ after a convenient choice of the determinations of argument and logarithm.
  Expand $\gvf(\xi-t)\psi(t)$ as a Taylor series $\displaystyle \sum_{m\geq 0} c_m(\xi) \, (t-\go)^m$ in $D_\go$. Since ${}_a\hspace{-.2em}\gc_\xi$ lies in $D_\go$ then, after commutation of sum and integral, $g_a(\xi)$ becomes 
\begin{equation*}
g_a(\xi)= \sum_{m\geq 0} c_m(\xi) \, \int_a^\xi (t-\go)^{m+\gl-1}\, \ln^p(t-\go)\, dt.
\end{equation*}
One can check that $\displaystyle \int_a^\xi (t-\go)^{m+\gl-1}\, \ln^p(t-\go)\, dt=G_{m+\gl,p}(\xi)-G_{m+\gl,p}(a)$ where
\begin{equation*}
G_{m+\gl,p}(\xi)=\left\lbrace \begin{array}{lll} \displaystyle
p!\,(\xi-\go)^{m+\gl}\, \sum_{s=0}^p\frac{(-1)^{p-s}}{s!} \,\frac{\ln^s(\xi-\go)}{(m+\gl)^{p-s+1}} \, & { \rm if }& m+\gl\neq 0\\
\noalign{\smallskip}
\displaystyle\frac{1}{p+1}\, \ln^{p+1}(t-\go)& {\rm if }& m=\gl=0\\
\end{array}\right.
\end{equation*}
For all $s=0,\dots,p$, the series $\displaystyle \sum_{m>0}c_m(\xi)\,\frac{(t-\go)^{m+\gl}}{(m+\gl)^{p-s+1}}$ converges on $D_\go$. Hence, the result.\qed

\begin{lemma}\label{antider} {\bf Anti-derivation}

Let $Q(X)\in\sM_{m,p}(\C\{\xi\}[X])$ be a polynomial matrix  with holomorphic coefficients at $\xi=0$.
We assume that the degrees of the successive columns of $Q$ are given by the row matrix $N_Q=[\nu\ \ \nu+1\ \cdots \ \nu+p-1]$.\\
Let $\Lambda\not\in-\N^*$.

The matrix function $\displaystyle \xi^\Lambda \xi^{-J_m}\,
Q(\ln \xi)\, \xi^{J_p}$ admits a unique anti-derivative of the form
\begin{equation*}
K(\xi)=\xi^{\Lambda+1} \xi^{-J_m}\, R(\ln \xi)\, \xi^{J_p}
\end{equation*}
where $R(X)\in\sM_{m,p}(\C\{\xi\}[X])$ is a polynomial matrix  with holomorphic coefficients at $\xi=0$ and same column-degrees $N_R=N_Q$ as $Q$.
\end{lemma}

Note that  the power of $\xi$ which can be factored increases by 1.
\medskip

\proof Denote $\displaystyle Q(X)=\sum_{k=0}^N Q_k\ln^k X$ and $\displaystyle R(X)=\sum_{k=0}^N R_k\ln^k X$ where $N=\nu+p-1$.

The derivative of $K$ reads 
\begin{equation*}
\begin{array}{ccl}
K'(\xi)&=&\xi^\Lambda \xi^{-J_m}\Big\lbrack
\big((\Lambda+1)I_m-J_m\big)(R_N\ln^N\xi+\cdots+R_1\ln\xi+R_0)
\\
\noalign{\medskip}
&&+\xi(R_N'\ln^N\xi+\cdots+R_1'\ln\xi+R_0') + (NR_N\ln^{N-1}\xi +\cdots+R_1)\\
\noalign{\medskip}
&&+(R_N\ln^N\xi+\cdots+R_1\ln\xi+R_0)J_{p}\Big\rbrack
\,\xi^{J_p}.
\end{array}
\end{equation*}
Identifying the powers of $\ln\xi$ we get the $N+1$ systems 
$$
\xi\, R_k'+\big((\Lambda+1)I_m -J_m\big) R_k+R_kJ_p=Q_k-(k+1)R_{k+1}\leqno{(*)_k} 
$$
{for}\ $k=0,1,\dots,N$ and $R_{N+1}=0$ which can be solved inductively starting with $k=N$. Like in Section \ref{sumressystem}, System $(*)_{k}$ can be solved iteratively  from the first to the last column and in each column, from the last row to the first one.
The fact that it admits a solution  holomorphic at 0 when $\gL\not\in-\N^*$ results from the fact that this is true for the differential equation $\xi \, y'+\gl\, y=a(\xi)$ when $a(\xi)$ is holomorphic at 0 and $\gl\not\in-\N$.\\
When $k$ is greater than $\nu$ set $k=\nu+\ell$. The fact that the first $\ell$ columns in $R_k$ can be chosen equal to 0 results from the fact that the same property holds for $Q_k$ and, by induction, for $R_{k+1}$. Then, it holds also for the right-hand side of System $(*)_k$ and the condition on the log-degree can be satisfied.

Unicity results from the fact that, since $\xi^{\Lambda+1}$ is neither a pole nor a constant, a non-zero constant cannot be put in such  form.
\qed

 \protect\boldmath
 \subsection{Singularities of $\gvf_m, m\geq 1$}\label{Singulphim} \protect\unboldmath
 Recall that the resurgent functions $\displaystyle \gvf_m=\left\lbrack \gvf_m^{j,\scriptscriptstyle\bullet}\right\rbrack_{1\leq j\leq J}$
 for $m\geq 1$ are iteratively determined, for all $j$, as  solutions of  the system 
 $$
(\xi-a_j)\frac{d\gvf_m^{j,\scriptscriptstyle\bullet}}{d\xi}-(\gl_j-1)\gvf_m^{j,\scriptscriptstyle\bullet}-J_{n_j} \gvf_m^{j,\scriptscriptstyle\bullet}+\gvf_m^{j,\scriptscriptstyle\bullet}J_{n_1} =\frac{d}{d\xi}\big(\hBjp{j}*\gvf_{m-1}\big)
\leqno{(\ref{systemphimblocs})}
 $$
 satisfying  convenient  initial conditions corresponding to those  satisfied by the $\tf_m$'s ($\gvf_0=\delta I_{n,n_1}$, $\gvf_m(0)=0$ for $m\geq  3$, $\gvf^{j,\scriptscriptstyle\bullet}_2(0)=0$ when $a_j\neq 0$ and a convenient non necessary 0 constant otherwise).
\medskip

Lemma \ref{antider} provides the log-degrees of the successive columns of an anti-derivative of the matrix $\xi^{\gl_k-1}\,\xi^{J_{n_k}}\, K\, \xi^{-J_{n_1}}$ when $K$ is a generic constant matrix, $0\leq {\rm Re}\,\gl_k<1$ and $\gl_k\neq 0$.
One can check that  the log-degrees are increased by 1 when $\gl_k =0$.
They are given  by 
\begin{equation}\label{N[k]}
N[k]=\left\lbrace
\begin{array}{ll}
\lbrack (n_k\hspace{-.2em}-\hspace{-.2em}1)\ \ (n_k\hspace{-.2em}-\hspace{-.2em}1)+1\;\cdots\;(n_k\hspace{-.2em}-\hspace{-.2em}1)+(n_1\hspace{-.2em}-\hspace{-.2em}1)\rbrack&\text{ if } \gl_k\neq 0\\
\noalign{\bigskip}
\lbrack \ n_k\quad\qquad n_k+1\qquad \cdots\qquad n_k+(n_1\hspace{-.2em}-\hspace{-.2em}1)\ \rbrack& \text{ if } 
\gl_k= 0.\\
\end{array}
\right.
\end{equation}

The behavior of $\gvf_m$ at a singular point  $\go\in\gO$ depends on the sheet of the Riemann surface $\sR_\gO$ we are on, \ie its depends on the path $\gc$ of analytic continuation followed from 0 (first sheet) towards $\go$. Following \cite{Sau05} we denote ${\rm cont}_\gc \gvf_m$ the analytic continuation of $\gvf_m$ along the path $\gc$.

\begin{dfn}\quad
{\rm \begin{itemize}
\item We say that $\xi\in\C\setminus \gO$ is {\it close  to} $\go\in\gO$ if there is a disc centered at $\go$ which contains $\xi$ and no other element of $\gO$ than $\go$.
\item We call {\it path from 0 towards} $\go\in\gO$ a path $\gc=\gc_\xi$ contained in $\C\setminus\gO$ which goes from 0 to a point $\xi$ close to $\go$.
\end{itemize}
}\end{dfn}

  \begin{pro} \label{singphim} Let $\go\in \gO$.
  
For any path of analytic continuation $\gamma$  from 0 towards $\go$, a  major of the singularity ${\buildrel {\scriptscriptstyle\,\, \nabla_\gc}\over {\gvf_m}}$ of ${\rm cont}_\gc\gvf_m$ at $\go$ exists  in the form
\begin{equation*}
\widecheck\gvf_m^{j;\scriptscriptstyle\bullet}(\go+\xi)= \xi^{\gl_j-1}\,\xi^{J_{n_j}}\, k_{m;(\go)}^{j;\scriptscriptstyle\bullet}\, \xi^{-J_{n_1}}+ {\rm rem}_m^{j;\scriptscriptstyle\bullet}(\xi)\quad \text{for all } j=1,\dots,J,
\end{equation*}
with a remainder 
$\displaystyle 
{\rm rem}_m^{j;\scriptscriptstyle\bullet}(\xi)=\sum_{\gl_k\vert a_k=\go} \xi^{\gl_k} \, R_{\gl_k,m}^ {j;\scriptscriptstyle\bullet}(\ln\xi)
$
 where 
\begin{itemize}
\item $k_{m;(\go)}^{j;\scriptscriptstyle\bullet}$ denotes a constant $n_j\times n_1$-matrix (recall that $n_k$ is the size of the $k^{th}$ Jordan block of the matrix $L$ of the exponents of formal monodromy) and 
 $k^{j;\scriptscriptstyle\bullet}_{m;(\go)}=0$ when $a_j\neq \go$,
\item  $R_{\gl_k,m}^{j;\scriptscriptstyle\bullet}(X)$ denotes a polynomial matrix with holomorphic coefficients at 0, the columns of which are of  degree $N[k]$ (\cf notation just above).
\end{itemize}
\end{pro}

	Of course,  $\widecheck\gvf_m^{j;\scriptscriptstyle\bullet}, k_{m;(\go)}^{j;\scriptscriptstyle\bullet}$
and ${\rm rem}_m^{j;\scriptscriptstyle\bullet}$ depend on  $\gc$ even though, for seek of simplicity, the notations  do not  show  it up.
\bigskip

Note that, in the remainders, the initial factor $\xi$ appears at powers like $\gl_k$ and no more $\gl_k-1$ so that the power $-1$ never occurs. Note also that, whatever are the values of $\gl_j$ and $\gl_k$, we have $\gl_j-1<\gl_k$ and this is why the terms in the remainders will always appear as subdominant.
\medskip

\proof  For all $m\geq 1$, System (\ref{systemphimblocs}) can be seen as an inhomogeneous linear system in the entries of $\gvf_m$. For any $j$, the general solution of the homogeneous system reads $(\xi-a_j)^{\gl_j-1} (\xi-a_j)^{J_{n_j}} \, k_{m;(\go)}^{j;\scriptscriptstyle\bullet}\, (\xi-a_j)^{-J_{n_1}}$ with $k_{m;(\go)}^{j;\scriptscriptstyle\bullet}$ an arbitrary constant $n_j\times n_1$ matrix and we have to prove that there is a particular solution of the inhomogeneous system in the form of the remainder. To this end, we integrate the system using the Lagrange method (variation of constants).  System (\ref{systemphimblocs}) is of the form $\displaystyle (\xi-a_j)\frac{d\phijp{j}_m }{d\xi} + {\rm linear \ terms } =Q$, the type of $Q$ depending upon $m$. Looking for a solution in the form $\gvf=(\xi-a_j)^{\gl_j-1} (\xi-a_j)^{J_{n_j}} \, K\, (\xi-a_j)^{-J_{n_1}}$ we obtain to determine $K$ up to a constant the condition $\displaystyle \frac{dK}{d\xi}=(\xi-a_j)^{-\gl_j} (\xi-a_j)^{-J_{n_j}} \, Q\, (\xi-a_j)^{J_{n_1}}$ and we just have to find anti-derivatives for the various possible $Q$.

For $m=1$, the inhomogenuity $\displaystyle Q=\frac{d\hBjp{j}}{d\xi}$ is an entire function, hence holomorphic at $\go$. When $a_j\neq \go$, then $\go$ is an ordinary point for the $j^{th}$ block of System (\ref{systemphimblocs}) and the inhomogenuity is holomorphic. Hence,  there is a holomorphic solution at $\go$ and we can choose $\widecheck\gvf_1^{j,\scriptscriptstyle\bullet}=0$. When $a_j=\go$,   Lemma \ref{antider} provides $K$ in the form $K=(\xi-a_j)^{-\gl_j+1} (\xi-a_j)^{-J_{n_j}} \, R(\xi)\, (\xi-a_j)^{J_{n_1}}$ with $R$ holomorphic at $\go$  and then, a particular solution $\phijp{j}_1(\xi)=R(\xi)$ holomorphic at $\go$ so that we can choose the remainder ${\rm rem}_1^j=0$. 

For $m=2$, using the superposition principle, we have to consider inhomogenuities $Q$ of the form  $\displaystyle Q=\frac{d \widehat B^{j;k}}{d\xi}*\phijp{k}_1  $ for all $k=1,\dots,J$. When $a_k\neq \go$, the solution $\phijp{k}_1  $ and then also $Q$ is holomorphic at $\go$ and we can conclude like in the case when $m=1$. From now, we forget about holomorphic terms. When $a_k= \go$, then 
$\phijp{k}_1  $ differs from a holomorphic function by terms of the form 
$(\xi-\go)^{\gl_k-1} (\xi-\go)^{J_{n_k}} \, k_{1;(\go)}^{k;\scriptscriptstyle\bullet}\, (\xi-\go)^{-J_{n_1}}$.
It results from Lemma \ref{convlog}   that, modulo a holomorphic function, $Q$ takes the form $(\xi-\go)^{\gl_k} P\big(\ln(\xi-\go)\big)$ with $P$ a $n_j\times n_1$-matrix of polynomials  with holomorphic coefficients the columns of which have log-degree $N[k]$. In both cases, $a_j=\go$ or $a_j\neq \go$, Lemma \ref{antider} provides a corresponding solution of the form $(\xi-\go)^{\gl_k} R\big(\ln (\xi-\go)\big)$ where $R$ is a polynomial matrix with holomorphic coefficients at $\go$ and column-log-degrees $N[k]$.

For $m\geq 3$, the inhomogenuity contains terms of the form $\displaystyle Q=\frac{d \widehat B^{j;k}}{d\xi}*\phijp{k}_{m-1}$.  
The factor $\phijp{k}_{m-1}$ splits into two parts: the first part coming from the general homogeneous solution is treated like in the case when $m=2$; the second part coming from the remainder is of a similar type but the fact that  all powers  $\lambda-1$ have been changed into $\lambda$. Consequently, there appears no factor $(\xi-\go)^{-1}$ and this insures the stability of the log-degree $N[k]$ since, from now, convolution generates no increase of the log-degree $N[k]$.
\qed

{\protect\boldmath}
 \subsection{Singularities of  {$\hf$}}
 {\protect\unboldmath}
 We are now ready to make explicit the form of the singularities of $\hf$ as a consequence of Theorem \ref{thmsum-res}, Proposition \ref{singphim} and an iterated application of the variation.
 
 Lemma \ref{lemmevar}  below states, without proof, some useful elementary properties of the variation. We denote ${\rm var}^p=\underbrace{{\rm var}\circ{\rm var}\circ\cdots\circ{\rm var}}_{p \text{ times}}$.

  \begin{lemma}\label{lemmevar}  \quad
  \begin{enumerate}
\item $\displaystyle {\rm var}\Big(\frac{\ln \xi}{2\pi i}\Big)=1$.
\item For all $p\in\N$,\\
$\begin{array}{ccl}
\displaystyle {\rm var}\Big(\Big(\frac{\ln \xi}{2\pi i}\Big)^p\Big)&=&\displaystyle \sum_{r=0}^{p-1}
(-1)^{p-r-1}C_p^r\Big(\frac{\ln \xi}{2\pi i}\Big)^r\\
\noalign{\medskip}
&=& \displaystyle p\Big(\frac{\ln \xi}{2\pi i}\Big)^{p-1} + \text{ lower log-degree terms}.\\
\noalign{\medskip}
\end{array}$\\
$\displaystyle \text{Consequently, }{\rm var}^p\Big(\Big(\frac{\ln \xi}{2\pi i}\Big)^p\Big)=p!\quad \text{and}\quad \displaystyle {\rm var}^{p+1}\Big(\Big(\frac{\ln \xi}{2\pi i}\Big)^p\Big)=0.$
\item For all $\gl\in\C$, $\displaystyle {\rm var}(\xi^\gl)=(1-e^{-2\pi i\gl})\xi^\gl$.\\
Consequently,  $\displaystyle {\rm var}^p(\xi^\gl)=(1-e^{-2\pi i\gl})^p\xi^\gl$ for all $p\in\N$ and \\ ${\rm var}(\xi^\gl)=0$ for all $\gl\in\Z$.
\item ${\rm var}(fg)={\rm var}(f) \, g+f\,{\rm var}(g) -{\rm var}(f)\,{\rm var}(g)$.\\
In particular, ${\rm var}(fg)=f\, {\rm var}(g)$ when ${\rm var}(f)=0$.
\item For all $\gl\in \C$ and $p\in\N$,
$$\displaystyle  \begin{array}{ccl}
\displaystyle {\rm var}\Big(\xi^\gl\Big(\frac{\ln \xi}{2\pi i}\Big)^p\Big)&=&\displaystyle (1-e^{-2\pi i\gl})\xi^\gl\Big(\frac{\ln \xi}{2\pi i}\Big)^p
+e^{-2\pi i\gl} \xi^\gl  {\rm var}\Big(\Big(\frac{\ln \xi}{2\pi i}\Big)^p\Big)\\
\noalign{\medskip}
&& \hspace{-2cm}\displaystyle =\ (1-e^{-2\pi i\gl})\xi^\gl\Big(\frac{\ln \xi}{2\pi i}\Big)^p
+ \xi^\gl \times \text{lower log-degree terms}.
\end{array}$$
Consequently, 
$$
{\rm var}^p\Big(\xi^\gl\Big(\frac{\ln \xi}{2\pi i}\Big)^p\Big)=(1-e^{-2\pi i\gl})^p\xi^\gl\Big(\frac{\ln \xi}{2\pi i}\Big)^p
+ \xi^\gl \times \text{lower log-degree terms}.
$$
\end{enumerate}
\end{lemma}

\begin{thm}\label{thmsingf} {\bf Singularity of \protect\boldmath $\hf$ at $\xi=\go$ \protect\unboldmath}\\
Let $\go\in \gO$.\\
 For any path of analytic continuation $\gc$ from 0 towards $\go$, a  major of the singularity ${\buildrel {\ \scriptscriptstyle\nabla_\gc}\over {f_\go}}$ of ${\rm cont}_\gc\hf$ at $\go$ exists  in the form
\begin{center} 
\noindent\begin{tabular}
{| p{8cm}  |}
\hline  
$\displaystyle \widecheck{f}^{j;\scriptscriptstyle\bullet}(\go+\xi)= \xi^{\gl_j-1}\,\xi^{J_{n_j}}\,k_{(\go)}^ {j;\scriptscriptstyle\bullet}\, \xi^{-J_{n_1}}+ {\rm Rem}_{(\go)}^{j;\scriptscriptstyle\bullet}(\xi)\tvi{18}{13}$ 
\\
  \hline
\end{tabular}\quad for all $j=1,\dots,J$
\end{center}
with a remainder 
$\displaystyle 
{\rm Rem}_{(\go)}^{j;\scriptscriptstyle\bullet}(\xi)=\sum_{\gl_\ell\vert a_\ell=\go} \xi^{\gl_\ell} \, R_{\gl_\ell;(\go)}^{j;\scriptscriptstyle\bullet}(\ln\xi)
$
 where 
\begin{description}
\item{$-$} \  $k_{(\go)}^{j;\scriptscriptstyle\bullet}$ denotes a constant $n_j\times n_1$-matrix (recall that $n_k$ is the size of the $k^{th}$ Jordan block of the matrix $L$ of the exponents of formal monodromy) and 
 $k_{(\go)}^{j;\scriptscriptstyle\bullet}=0$ when $a_j\neq \go$,
\item{$-$} \   $R_{\gl_\ell;(\go)}^{j;\scriptscriptstyle\bullet}(X)$ denotes a polynomial matrix with {\rm summable-resurgent coefficients} in $\widehat{Res}_{\gO-\go}^{\rm sum}$, the columns of which are of  log-degree $N[\ell]$ (\cf Section \ref{Singulphim}, Formula (\ref{N[k]})).
\end{description}
\end{thm}
Recall that the major $\displaystyle \widecheck{f}^{j;\scriptscriptstyle\bullet}$,  the constant matrix
$k_{(\go)}^{j;\scriptscriptstyle\bullet}$ and  the coefficients of the remainder 
$R_{\gl_\ell;(\go)}^{j;\scriptscriptstyle\bullet}$ depend on  $\gc$ even though, for seek of simplicity, we do not show it up in the notations. Note that, in practice, $k_{(\go)}^{j;\scriptscriptstyle\bullet}$ can be determined as the coefficient of the monomial $\xi^{\gl_j-1}$.
\bigskip

It is worth to make explicit the following two particular cases.
\begin{itemize}
\item {\bf Case with diagonal formal monodromy: \protect\boldmath $\displaystyle L=\oplus_{j=1}^n\gl_j$\protect\unboldmath }

In this case, $\widecheck{f}^{j;\scriptscriptstyle\bullet} $ reduces to just one entry which we denote        $\widecheck{f}^{j}$.
\begin{center} 
\noindent\begin{tabular}
{| p{10cm}  |}
\hline  
$\displaystyle \widecheck{f}^{j}(\go+\xi)= 
\frac{k_{(\go)}^j}{\xi}+h_{(\go)}^j(\xi)\,\ln\xi+ \sum_{\gl_\ell\neq 0\,\vert\,a_\ell=\go}
\xi^{\gl_\ell-1}\, h_{\gl_\ell;(\go)}^j(\xi)\tvi{18}{13}$
\\
  \hline
\end{tabular}
\end{center}
where $k_{(\go)}^j$ is a constant, $h_{(\go)}^j$ and all $h^j_{\gl_\ell;(\go)}$ are summable-resurgent functions of $\widehat{Res}_{\gO-\go}^{\rm sum}$. Moreover, 

\begin{itemize}
\item \  $k_{(\go)}^j=0 $  when  $a_j\neq\go$  or  $\gl_j\neq 0$,
\item \  $h_{(\go)}^j=0$  when there exists no  $\ell$  such that 
 $a_\ell=\go$  and  $\gl_\ell=0$.
 \end{itemize}

 \renewcommand{\arraystretch}{1} 

\item {\bf Case with trivial formal monodromy: \protect\boldmath $\displaystyle L=O_n$\protect\unboldmath }
\begin{center} 
\noindent\begin{tabular}
{| p{5.5cm}  |}
\hline  
$\displaystyle \widecheck{f}^{j}(\go+\xi)= 
\frac{k_{(\go)}^j}{\xi}+h_{(\go)}^j(\xi)\,\ln\xi\tvi{18}{13}$
\\
  \hline
\end{tabular}
\end{center}
where $k_{(\go)}^j$ is a constant and $h_{(\go)}^j$ a summable-resurgent function of $\widehat{Res}_{\gO-\go}^{\rm sum}$. Moreover, 

\begin{itemize}
\item \  $k_{(\go)}^j=0 $  when  $a_j\neq\go$,
\item \  $h_{(\go)}^j=0$  when there exists no  $\ell$  such that 
 $a_\ell=\go$.
 \end{itemize}
 \end{itemize}
 
Thus, under the condition that there exists no $\gl_\ell\neq 0$ associated with  $a_\ell=\go$, the singularity ${\buildrel {\scriptscriptstyle\nabla}\over {f_\go}}$ of $\hf$ at $\go$ is simple in the restrictive sense (\cf Definition \ref{footnotesimple} and its footnote).
\vspace{.5cm}

\proof --- Prove first the weaker result that asserts that the singularity of $\hf$ at $\go$ has the given form  with holomorphic coefficients at 0. \\
Given $0<\nu<\nu_1$ and $\gve>0$, we consider the domains
\begin{center}
$\gD_\nu=\{\nu<\vert \xi\vert<\nu_1\}$ in $\C$\tvi{0}{15} \\
and\quad $\gD_\nu'=\{\nu<\vert\xi\vert<\nu_1, \ \gt-\gve-2\pi N'<{\rm arg}(\xi)<\gt+\gve +2\pi\}$
\end{center}
a lift of the ring $\gD_\nu$ to $N'+1$ consecutive sheets of the Riemann surface of the logarithm at 0. The number $N'$ will be determined later; it has to be finite and large enough. The argument $\gt$ fixes the sheet on which $\hf$ is studied.\\
We assume $\nu_1$ so small that the disc $\vert \xi\vert\leq \nu_1$ lies at a distance at least $\nu$ from $(\gO-\go)\setminus \{0\}$.

We fix $j$. From Proposition \ref{singphim}  we can write 
$$
{\rm cont}_\gc \gvf_m^{j,{\scriptscriptstyle\bullet}}(\go+\xi)=\sum_{\rm \gl\not\in\Z}\sum_{p=0}^{p_\gl}h_m^{\gl,p}(\xi)\,\xi^\gl\,\ln^p\xi+\frac{1}{\xi}\sum_{p=0}^{p_0} h_m^{0,p}(\xi)\,\ln^p\xi
$$
where the pairs $(\gl,p)$ are supposed distinct modulo $\Z$ and finitely many.\\
Denote $\gL$ the set of all exponents $\gl$ appearing in these summations. Elements in $\gL$ are equal to either a $\gl_\ell$ or $\gl_j-1$ when $\gl_j\not \in\Z$ and $a_j=\go$. When $\gl_j=0$ and $a_j=\go$, instead of introducing $\gl=\gl_j-1=-1$ we factor $\displaystyle\frac{1}{\xi}$ so that no polar part occurs in the coefficients $h_m^{0,p}(\xi)$.  We denote $p_\gl$ the highest log-degree $p$ associated with $\xi^\gl$.\\ 
From Proposition \ref{singphim} we know that all coefficients $h_m^{\gl,p}(\xi)$ are holomorphic at 0 and we have to prove that $\hf(\go+\xi)$  has the same form as the $\gvf_m$'s  and holomorphic coefficients at 0. The proof mostly relies on the uniform convergence of the series $\sum_{m\geq 1} \gvf_m$ to $\hf-\gd I_{n,n_1}$ on any compact set avoiding $\gO$ (\cf Remark \ref{unifconv}). For simplicity, we skip writing ``${\rm cont}_\gc$'' although we consider analytic continuations along $\gc$.
\bigskip

Suppose, to begin, that all majors $\widecheck\gvf_m^{j,\scriptscriptstyle\bullet}$ are equal to 0 which means that all  $\displaystyle \gvf_m^{j,\scriptscriptstyle\bullet}(\go+\xi)$ are  holomorphic functions $h_m^{j,\scriptscriptstyle\bullet}(\xi)$ at 0 and then, all holomorphic on the same disc $\vert\xi\vert<\nu_1$ previously chosen.
Since the series $\sum_{m\geq 0}\gvf_m^{j,\scriptscriptstyle\bullet}(\go+\xi)=\sum_{m\geq 0} h_m^{j,\scriptscriptstyle\bullet}(\xi)$ converges uniformly to $\hf^{j,\scriptscriptstyle\bullet}(\go+\xi)$ on compact sets avoiding $\gO-\go$ (Theorem \ref{thmsum-res}, Remark \ref{unifconv}) 
the function $\hf^{j,\scriptscriptstyle\bullet}(\go+\xi)=\sum_{m\geq 0} h_m^{j,\scriptscriptstyle\bullet}(\xi)$ is a holomorphic function on the punctured disc $0<\vert \xi\vert<\nu_1$.
  Its Laurent series is the sum of the Laurent series of the $ h_m^{j,\scriptscriptstyle\bullet}(\xi)$'s; hence, it displays no polar part and $\hf^{j,\scriptscriptstyle\bullet}(\go+\xi)$  can be continued into a holomorphic function $h^{j,\scriptscriptstyle\bullet}(\xi)$  on the disc $\vert \xi\vert<\nu_1$. We can conclude, in this case, that a major of $\hf$ is also equal to zero: $\widecheck{f}^{j,\scriptscriptstyle\bullet}(\go+\xi)=0$.

Given $(\gl,p)$, to prove that the series $\sum_{m\geq 0}h_m^{\gl,p}(\xi)$ converges to a holomorphic function $h^{\gl,p}(\xi)$ about $\go$ we proceed as above after having reduced all $\gvf_m^{j,\scriptscriptstyle\bullet}(\go+\xi)$ to $h_m^{\gl,p}(\xi)$ by means of iterated variations as indicated below.
 We  can  then  conclude by the same arguments as above since the uniform convergence property keeps valid for the variations as well.

We base the reduction on the properties of the variation stated in Lemma \ref{lemmevar}:  the variation of $\ln^p\xi$ iterated $p$ times produces  a non-zero constant and it produces 0 in one more step; the iterated variation of $\xi^\gl\,\ln^p\xi$ when $\gl\not\in\Z$ generates a dominant term of the same form times a non zero constant. 

Fix $\gl\in\gL$. \\
 It is sufficient to consider the case of the monomials $h_m^{\gl,p_\gl}(\xi)\,\xi^\gl\,\ln^{p_\gl}\xi$ of highest log-degree. Indeed, we can then proceed iteratively  on the descending log-degrees terms after cancellation  of the highest log-degree terms.\\
 Here is a possible way to reduce $\gvf_m^{j,\scriptscriptstyle\bullet}(\go+\xi)$ to $h_m^{\gl,p_\gl}(\xi)$:\\
 For each $\gl'\in\gL$, $\gl'\neq\gl$ successively, multiply by $\xi^{-\gl'}$, take the variation $p_{\gl'}+1$ times and multiply by $\xi^{\gl'}$, thus canceling all terms factored by $\xi^{\gl'}$. We are left with only terms in $\xi^\gl \, \ln^p\xi$, for $p=0,1,\dots,p_\gl$. Multiply by $\xi^{-\gl}$ and take the variation $p_\gl$ times. We get so $h_m^{\gl,p_\gl}(\xi)$ up to a non-zero constant.\\
 We can now estimate a convenient value for $N'$: the process, to be valid, requires that $N'$ be as large as the total number of variations used.

 This ends the proof of the fact that the singularity of ${\rm cont}_\gc\hf$ at $\go$ has the given form  with holomorphic coefficients at 0. In particular, the constant matrices $k_{(\go)}^{j;\scriptscriptstyle \bullet}$ are given by $k_{(\go)}^{j;\scriptscriptstyle \bullet}=\sum_{m\geq 1}k_{m:(\go)}^{j;\scriptscriptstyle \bullet}$.
 \bigskip
 
 --- The proof of the fact that the coefficients $h_{(\go)}^{\gl,p}(\xi)$ are actually  summa\-ble-resurgent functions in $\widehat{\sR es}_{\gO-\go}^{\rm sum}$ is straightforward from the fact that their germs at the origin are equal to iterated variations of functions themselves in $\widehat{\sR es}_{\gO-\go}^{\rm sum}$. \qed
 
 \subsection{Principal major and connection constants}
Let $\gt\in\R/2\pi\Z$ be  an anti-Stokes direction and  $\go\in\gO_\gt$ a Stokes value in direction $\gt$ associated with $\tf(x)$.

The constants $k_{(\go)}^{j;\scriptscriptstyle{\bullet}}$ and the polynomials $R_{\gl_\ell,(\go)}^{j;\scriptscriptstyle{\bullet}}$ found in Theorem \ref{thmsingf} depend, as already said, on the path of analytic continuation $\gc_\xi$ and meanwhile, on the chosen determination  of the argument  around $\go$. \\

$\bullet$ We consider a path $\gc^+$ from 0 towards $\go$ defined as follows:

 \noindent\begin{minipage}[c]{85mm}
\baselineskip=1.5em
 $\gc^+=\gc^+_\xi \tvi{15}{0}$ goes   along the straight line $[ 0,\go]$ from 0 towards $\go$ and bypasses all  intermediate singular points $\go'\in\gO_\gt\cap ]0,\go ]$ to the right as shown on the figure. 
\end{minipage}
\hskip 1cm
\begin{minipage}[c]{60mm}
\includegraphics[scale=1]{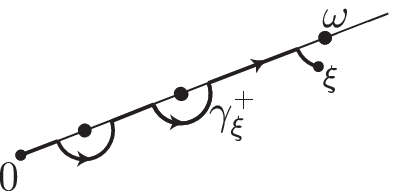}\tvi{0}{10}
\end{minipage}\\

$\bullet$ We choose the principal determination of the variable $\xi$ around $\go$\footnote{\label{determ} Any choice is convenient. However, to be compatible, on the Riemann sphere, with the usual choice $0\leq \arg(\gz=1/\xi)<2\pi$ of the principal determination  at infinity we suggest to choose $-2\pi<\arg(\xi) \leq 0$ as principal determination about 0 as well as about any $\go$ at finite distance.}.
\medskip

  Theorem \ref{thmsingf} allows to set the following definition.

\begin{dfn} {\bf Principal major}\\
{\rm We call {\it principal major of} $\hf$ {\rm at } $\go\in\gO_\gt$ the major $\widecheck f^+(\go+\xi)$  of ${\rm cont}_{\gc^+}(\hf)$  when the principal determination of the variable $\xi$ is chosen.
\\
According to Theorem \ref{thmsingf}, for  $j=1,\dots,J$ and convenient constant matrices  $k_{(\ugo)}^ {+^{j;\scriptscriptstyle\bullet}}$ it reads
\begin{equation}
 \widecheck{f}^{+^{j;\scriptscriptstyle\bullet}}(\go+\xi)= \xi^{\gl_j-1}\,\xi^{J_{n_j}}\,k_{(\ugo)}^ {+^{j;\scriptscriptstyle\bullet}}\, \xi^{-J_{n_1}}+\sum_{\gl_\ell\vert a_\ell=\go} \xi^{\gl_\ell} \, R_{\gl_\ell;(\ugo)}^{+^{j;\scriptscriptstyle\bullet}}(\ln\xi)\tvi{18}{13}
\end{equation}
 Moreover, in case there exists $\ell$ such that $\gl_\ell=0$ and $a_\ell=\go$, we assume that
$R_{0;(\ugo)}^{+^{j;\scriptscriptstyle\bullet}}(0) =0$ so that the principal major never contains a holomorphic term.}
\end{dfn}

With the requested  conditions the principal major of $\hf$ at $\go$ is uniquely determined. We denote  $k_{(\ugo)}^ {+}$  with a dot to indicate the choice of the principal determination.

By means of elementary transformations (\cf Section \ref{prepared}) we can extend the previous results and definitions to any of the $k^{th}$ column-blocks of $\hF(\xi)$ changing the Stokes values $a_\ell$   into $a_\ell-a_k$ (actually, with respect to $\hf(\xi)$, $a_\ell$ stands for $a_\ell-a_1$ with $a_1=0$)   and the exponents of formal monodromy $\gl_\ell$ into $\gl_\ell-\gl_k$. We can then reformulate Theorem \ref{thmsingf} as follows.
\begin{cor} Let $\gt\in\R/2\pi \Z$ be an anti-Stokes direction and  $\go\in{\bf \gO}_\gt$  a Stokes value in direction $\gt$ associated with System (\ref{A}).

To the choice of the principal determination of the variable $\xi$ there is a unique constant matrix $K^+_{(\ugo)}$ such that the principal major $\widecheck F^+$ of $\hF$ at $\go$ reads  
\begin{equation} \label{ppalmajor}
\renewcommand{\arraystretch}{2} 
\begin{tabular}{| p{8cm}|} 
\hline 
$\quad\displaystyle  \widecheck F^+(\go+\xi)= \frac{1}{\xi}\,\xi^L\, K^+_{(\ugo)}\,\xi^{-L} + {\rm Rem}^+(\ln\xi)
\tvi{20}{18}
$\\
\hline
\end{tabular}\,
\end{equation}
The matrix $ K^+_{(\ugo)}$ is a constant $n\times n$-matrix satisfying  $K^{+^{j;k}}_{(\ugo)}=0$ for all $(j;k)$ such that  $a_j-a_k\neq\go$. In particular, all diagonal blocks are equal to zero.\\
The remainder ${\rm Rem}^+(\ln\xi)$ is a linear combination of polynomials in $\ln\xi$ with summable-resurgent coefficients weighted by subdominant complex powers of $\xi$. 
\end{cor}

\begin{dfn} {\bf Matrix of the connection constants}\\
{\rm The $n\times n$-matrix }
\begin{equation} \label{connconst}
\renewcommand{\arraystretch}{2} 
\begin{tabular}{| p{3.5cm}|}
\hline 
$\quad\displaystyle  K_{\ugt}^+=\sum_{\go\in{\bf \gO}_\gt} K^+_{(\ugo)}
\tvi{20}{20}
$\\
\hline
\end{tabular}
\end{equation}
{\rm is called} {\it (principal) matrix of the connection constants of} 
$\hF$ {\rm in direction} $\gt$. 
\end{dfn}

\noindent A $n_j\times n_k$-block $K_{\ugt}^{+^{j;k}}$ of $K_{\ugt}^+$ is equal to 0 when $a_j-a_k$ does not belong to ${\bf \gO}_\gt$.
The possibly non-zero entries of $K^+_{\ugt}$ are also called {\it principal  multipliers of connection in direction $\gt$}.
\bigskip

Note that we still need no other structure on $\widehat{\sR es}_\gO$ and $\sC$ than the structure of $\sO^{\leq 1}(\C)*$module and, in particular, there is no need yet to develop a structure of a convolution algebra.


\section[Merom. inv.: Stokes cocycle or alien derivatives]{Meromorphic invariants: Stokes cocycle or alien derivatives}\label{MeroInv}

The classifying set for the local meromorphic classification at 0 of connections endowed with an isomorphism of their formalized (\ie  formal gauge transformations $\hF$) in a given formal class was given by Y. Sibuya \cite{Sib77, Sib90} and B. Malgrange \cite{Mal83} in terms of a non-Abelian 1-cohomology set. Actually, this classifying set can be given  a structure of a unipotent Lie group and it is  isomorphic to the direct product of the Stokes-Ramis groups in each anti-Stokes direction associated with the connection (\cf  D.G. Babbitt and V.S. Varadarajan \cite{BV89} for an abstract proof and M. Loday-Richaud \cite{L-R94} for a constructive one). Indeed, in each 1-cohomology class there exists a unique special cocycle, called {\it Stokes cocycle}, whose components coincide with the Stokes-Ramis automorphisms independently defined as the defects of analyticity of $\hF$ in each anti-Stokes direction \cite{R85};  \cite[ Prop. III.2.1, Th. III.2.8]{L-R94}.  

To any formal class and associated anti-Stokes direction the Stokes-Ramis automorphisms form a free Lie group conjugate to a group of unipotent triangular matrices submitted to some vanishing conditions (\cite{L-R94} Def. I.4.12). Its Lie algebra is conjugate to an algebra of nilpotent matrices submitted to similar vanishing conditions and the exponential map sends it homeomorphically onto the Stokes-Ramis group.
 It is then equivalent to characterize a  meromorphic class by giving its image in the Stokes-Ramis groups (\ie{\it its Stokes cocycle} or {\it its Stokes matrices} after the choice of a $\C$-basis of solutions) or its image in the Lie algebra, tangent space of the Stokes-Ramis groups at the identity (\ie {\it its alien derivatives}).
\bigskip

Here below, we perform 
these descriptions in more details and we provide an explicit formula for the  Stokes matrices in the Laplace plane  in terms of the connection matrices in the Borel plane.
\bigskip

\noindent{\bf Note.} From now, each time a determination of the argument is required, we choose the principal determination fixed similarly in the Laplace and in the Borel plane (\cf Footnote \ref{determ}). Given an anti-Stokes direction $\gt\in\R/2\pi\Z$ we denote
$\gt^\star\in\R$ the chosen determination of $\gt$ and $\go^\star$ the Stokes value $\go$ in direction $\gt$ with the same determination. We keep denoting the variables $x, \xi,\dots$ while indicating $\arg(x)\simeq \gt^\star, \arg(\xi)\simeq \gt^\star, \dots$.

\subsection{Stokes cocycle}

Let us start with a description of the Stokes-Ramis automorphisms, components of the Stokes cocycle, from the viewpoint of summation.

Recall (\cf Section \ref{knownresults}) that the anti-Stokes directions attached to the first $n_1$ columns $\tf$ of the gauge transformation $\tF$ are defined as the directions of the various non-zero Stokes values  of System (\ref{A}), \ie the Stokes values belonging to $\gO^*=\{a_1,\dots,a_J\}\setminus \{0\}$ (see notations of Section 
\ref{knownresults}).
To the $\ell^{th}$ column-block of $\tF$ the set $\gO$ must be translated to $\gO-a_\ell$. We keep denoting ${\bf \gO}=\{a_j-a_\ell\neq 0\}$ the set of all non-zero Stokes values attached to $\tF$.\\ \\
Given $\gt\in\R/2\pi\Z$ let $d_\gt$ denote the half line issuing from 0 with argument $\gt$ and set $\gO_\gt=\gO^*\cap d_\gt$ and ${\bf \gO}_\gt= {\bf \gO}^*\cap d_\gt$.

\subsubsection{Stokes automorphisms as  gauge transformations}

When $\gt$ is not an anti-Stokes direction for $\tF$ (\ie ${\bf \gO}_\gt=\emptyset$) then $\tF$ can be applied a Borel-Laplace integral $\displaystyle \int_{d_\gt} \hF(\xi)e^{-\xi/x}d\xi$  in direction $\gt$ (\cf Theorem \ref{thmsum-res}, for instance) and in neighboring, not anti-Stokes, directions giving thus rise to an analytic function $s_\gt(\tF)$ defined and 1-Gevrey asymptotic to $\tF$ on a sector  $\gS_{\gt,>\pi}$ bisected by $\gt$ with opening larger than $\pi$. The function $s_\gt(\tF)$ is called {\it 1-sum} or {\it  Borel-Laplace sum} of $\tF$ in direction $\gt$.

When $\gt$ is an anti-Stokes direction for $\tF$ the Borel-Laplace integral does not exist anymore in general. However, taking the limit as  $\gve$ tends to 0 of the Borel-Laplace sums in directions $\gt-\gve$ and $\gt+\gve$ one defines, by analytic continuation, two analytic functions, 1-Gevrey asymptotic to $\tF$ on a germ of half-plane $\gS_{\gt,\pi}$ bisected by $\gt$. 
\smallskip

\noindent\begin{minipage}[c]{35mm}

\hspace{.7cm}
\includegraphics[scale=1]{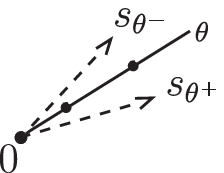}

\end{minipage}
\hspace{.2cm}
\begin{minipage}[c]{90mm}
\baselineskip=1.5em
We call {\it sum of $\tF$ to the right} of $\gt$ and we denote $s_{\gt^+}(\tF)$ the analytic continuation to $\gS_{\gt,\pi}$ of $s_{\gt-\gve}(\tF)$ as $\gve$ tends to 0. 
We call {\it sum of $\tF$ to the left} of $\gt$ denoted $s_{\gt^-}(\tF)$ the analytic continuation of $s_{\gt+\gve}(\tF)$.
\end{minipage}
\footnote{ The denominations  {\it  to the right} and {\it  to the left} fit the natural orientation around 0 on the Riemann sphere. Our choice of the signs $+$ and $-$ in $s_{\gt^+}$ and $s_{\gt^-}$ may look inappropriate to such an orientation but our will is to fit the usual notations at infinity. Indeed, positioning the singularity at 0 or at infinity exchanges the orientation on the Riemann sphere. The Borel transform at 0 does not exchange the orientation while it does it at infinity. In both cases, we can refer to the fact that the sums $s_{\gt^+}$ and $s_{\gt^-}$ correspond to Laplace integrals in the Borel plane along a path following a line $d_\gt$ and passing the Stokes values positively for $s_{\gt^+}$ and negatively for $s_{\gt^-}$ (\cf Fig. Section \ref{StoConnection}). } 
\\  

The {\it Stokes phenomenon} stems from the fact that the
 two lateral sums $s_{\gt^+}(\tF)$ and $s_{\gt^-}(\tF)$ are not analytic continuations from each other in general. The defect of analyticity is quantified by the Stokes automorphisms 
\begin{equation} 
\renewcommand{\arraystretch}{2} 
\begin{tabular}{| p{5cm}|}
\hline 
$\displaystyle
St_{\gt;\tF}=s_{\gt^-}(\tF)^{-1}\circ s_{\gt^+}(\tF)
$\\
\hline
\end{tabular}\,. 
\end{equation} 
in each anti-Stokes direction $\theta$.
Thus defined, the Stokes automorphism of $\tF$ in direction $\theta$ is an automorphism of the normal form (\ref{A0}); precisely, a gauge transformation, 1-Gevrey flat on $\Sigma_{\theta,\pi}$, which leaves invariant System (\ref{A0}). The Stokes automorphism $St_{\gt;\tF}$  depends on $\tF$ and $\gt$ and it does not depend on the choice of a determination of the argument near $\theta\in\R/2\pi\Z$.
\subsubsection{Matrix representations}

\hspace{.8cm} One can give the Stokes automorphism $St_{\gt;\tF}$  matrix representations in  $GL(n,\C)$ by associating with the (formal) normal solution $\displaystyle \tY_0(x)=x^Le^{Q(1/x)}$ an actual one near $\theta$ by means of a choice of a determination of the argument in the Laplace plane:  

We change the formal power $\displaystyle x^L$ into an actual function near $\theta$ by choosing the principal  determination of the argument. We get thus an actual function, still denoted $x^L$, defined and analytic for $\arg(x)$ close to $\ugt$ (denote $\arg(x)\simeq \ugt$).
 In our case of a single level equal to one, the polynomials $q_j(1/x)$ in $Q(1/x)$ do not require any choice of a determination of the argument. However, a formal exponential determines an actual exponential only up to a multplicative constant\footnote{\label{exptorus} What is called  formal exponential $\displaystyle e^p$ is ``the'' formal solution of the differential equation $y'-p'y=0$, which is well defined only up to a multiplicative scalar.}. We fix such a constant once for all, for instance, by choosing the function denoted $\displaystyle  e ^{Q(1/x)}$ in the usual analytic meaning. 
 We denote $Y_{0,\ugt}(x)$ the fundamental solution such defined for $\arg(x)\simeq \ugt$.
 
To the actual normal solution  $Y_{0,{\ugt}}(x)$, there cor\-res\-pond  two analytic fundamental solutions  $s_{\gt^-}(\tF)Y_{0,{\ugt}}(x)$ and $s_{\gt^+}(\tF)Y_{0,{\ugt}}(x)$  of System (\ref{A}).
 There exists then a unique constant invertible matrix (called {\it Stokes-Ramis matrix} or simply {\it Stokes matrix}\footnote{ In the literature, a {\it Stokes matrix} has often a more general meaning where one allows to compare any two asymptotic solutions whose domains of definition overlap. We exclude such a generality here.}), which we denote $I_n+C_{\ugt}$, such that 
 \begin{equation}\label{sum+-}
\displaystyle \,  s_{\gt^+}(\tF)\,Y_{0,{\ugt}}(x)=s_{\gt^-}(\tF)\,Y_{0,{\ugt}}(x)(I_n+C_{\ugt})\tvi{14}{8}
\end{equation}
 and the Stokes automorphism $St_{\gt;\tF}$ is given a matrix representation by the formula
\begin{equation} \label{Stokesauto}
\renewcommand{\arraystretch}{2} 
\begin{tabular}{| p{11cm}|}
\hline 
$\displaystyle \begin{array}{l}
\ St_{\gt;\tF}\equiv s_{\gt^-}(\tF)^{-1}\,s_{\gt^+}(\tF)=x^Le^{Q(1/x)}(I_n+C_{\ugt})e^{-Q(1/x)}x^{-L}\tvi{20}{0}\\
\hfill \text{for } \arg(x)\simeq \ugt\qquad\\
\end{array}$\\
\hline
\end{tabular}\,. 
\end{equation}

Thus, although the Stokes automorphisms, components of the Stokes cocycle, are intrinsically determined from $\tF$, their matrix representations as constant invertible matrices are defined only up to 
the choice of an actual normal solution $Y_{0,\ugt}(x)$ in each anti-Stokes direction. As it is clear from Formula (\ref{Stokesauto}) this indetermination results in the conjugacy action of the formal monodromy (change of ``formal determination'' $x\mapsto xe^{2\pi i}$ in the normal solution)  and of the exponential torus (see Section \ref{aliender}). Temporarily, we forget about the action of the exponential torus by systematically associating with a formal exponential $\displaystyle e^q$ the analytic function denoted the same way (choice of the arbitrary constant equal to 1) so that, a (formal) normal solution being given, our unique freedom lies in the choice of the determination of the argument of the variable $x$ in a neighborhood of the direction $\gt$, \ie the choice of a determination $\ugt$ of $\gt$.

The Stokes automorphism $St_{\gt;\tF}\tvi{15}{0}$ is unipotent (\ie its matrix $C_{\ugt}$ is nilpotent whatever the choice of a determination of $\ugt$) due to the fact that $s_{\gt^-}(\tF)$ and $s_{\gt^+}(\tF)$ are both 1-Gevrey asymptotic to the same matrix $\tF$. Indeed, a $n_j\times n_\ell$ block $C_{\ugt}^{j,\ell}\tvi{14}{0}$ is 0 as soon as $a_j-a_\ell\not\in\gO_\gt$. In particular, the diagonal blocks are equal to 0 and one can put $C_{\ugt} $ in triangular form by conveniently reordering the Stokes values $a_j$ (\cf \cite[Consequence I.4.8]{L-R94}).\\

  Formula (\ref{Stokesauto}) has an  ``additive'' form  
\begin{equation}\label{Stokesautoadd}
\renewcommand{\arraystretch}{1} 
\begin{tabular}{| p{11.5cm}|}
\hline 
$\displaystyle 
\begin{array}{l}
\ s_{\gt^+}(\tF)(x)-s_{\gt^-}(\tF)(x)=s_{\gt^-}(\tF)(x) \, x^L\,e^{Q(1/x)}\, C_{\ugt}\,e^{-Q(1/x)}\,x^{-L}
\tvi{20}{0}\\
\hfill \text{for } \arg(x)\simeq \ugt\qquad\\
\end{array} 
$\\
\hline
\end{tabular}
\end{equation}
which will be used later.

\subsubsection{Stokes automorphims acting on formal solutions}\label{formalsto}

The Stokes-Ramis matrix $I_n+C_{\ugt}$ was introduced in Formula (\ref{sum+-}) as the matrix of a linear  map
\begin{equation*}
s_{\gt^-}(\tF)\,Y_{0,\ugt}(x) \mapsto s_{\gt^+}(\tF)\,Y_{0,\ugt}(x)
\end{equation*}
associating sums to the right with sums to the left in direction $\ugt$. The map is defined and bijective in the space of actual solutions of System (\ref{A}) over a germ of half plane $\gS_{\gt,\pi}$.

Coming back from the actual normal solution $Y_{0,\ugt}(x)\tvi{15}{0}$ to the formal one $\displaystyle \tY_0(x)=x^Le^{Q(x)}$ the previous map can be read as a bijective linear map 
\begin{equation}\label{stoformel}
St_{\ugt} : \quad \tF(x)\, x^L\, e^{Q(1/x)} \mapsto \tF(x)\, x^L\, e^{Q(1/x)} (I_n+C_{\ugt})
\end{equation}
in the space of formal solutions of System (\ref{A}). Such a map depends on $\ugt$ and no more only  on $\gt$.

One also calls the map $St_{\ugt}$ a {\it Stokes automorphism} in direction $\ugt$. The Stokes cocycle of System (\ref{A}) is equivalent to the collection of the Stokes maps $St_{\ugt}$ for $\ugt$ running over  arguments of the anti-Stokes directions in a fundamental domain, say the principal one $-2\pi< \arg \leq 0$.
\bigskip

Note that, to the first column-block $\tf$ made of formal power series (recall $a_1=\gl_1=0$), there does not correspond a power series in general.  In Section \ref{aliender}, in order to define alien derivations, we will   extend such a map into an automorphism of a differential algebra containing the formal solutions of System (\ref{A}).


\protect\boldmath
\subsection{Maps $\dDodp$}\label{MapDelta+}
\protect\unboldmath

Given $\go\in{\bf \gO^*}$ recall that $\ugo$ denotes its principal determination  with argument $\ugt$.

The Stokes automorphism $St_{\ugt}$  can be split into a sum of linear maps  $\dDodp$, each of them taking into account the contribution of a different Stokes value $\go\in{\bf \gO}^*$ as follows: their matrices, also denoted $\dDodp$, are obtained from $C_{\ugt}$ by keeping unchanged the blocks $C_{\ugt}^{j,\ell}$ such that $a_j-a_\ell= \go$ and equating all other blocks to 0.\\
Obviously, $\displaystyle C_{\ugt} = \sum_{\go\in{\bf\gO}_\gt} \dDodp$ and Formula (\ref{Stokesautoadd}) reads,  for $\arg(x)\simeq \ugt$,
\begin{equation}
\renewcommand{\arraystretch}{2} 
\begin{tabular}{| p{10.8cm}|}
\hline 
$\displaystyle \ 
s_{\gt^+}(\tF)(x)-s_{\gt^-}(\tF)(x)= s_{\gt^-}(\tF)(x)\, \sum_{\go\in{\bf\gO}_\gt}\, e^{-\go/x}\,x^{L}\, \dDodp\, x^{-L}\tvi{0}{20}
$\\
\hline
\end{tabular}
\end{equation}

In restriction to the first $n_1$ columns and denoting  $\ddodp$ the first $n_1$ columns of $\dDodp$, we obtain, for $\arg(x)\simeq \ugt$,
\begin{equation}\label{stomatbis}
\renewcommand{\arraystretch}{2} 
\begin{tabular}{| p{10.6cm}|}
\hline 
$\displaystyle \ 
s_{\gt^+}(\tf)(x)-s_{\gt^-}(\tf)(x)= s_{\gt^-}(\tF)(x)\, \sum_{\go\in{\bf\gO}_\gt}\, e^{-\go/x}\,x^{L}\, \ddodp\, x^{-J_{1}}\tvi{0}{20}
$\\
\hline
\end{tabular}\,.
\end{equation}

\subsection{Stokes-Ramis  versus connection matrices} \label{StoConnection}

The left hand side of Formula (\ref{stomatbis}) can be seen as the Laplace integral 
\begin{equation}\label{stomatadd}
s_{\gt^+}(\tf)(x)-s_{\gt^-}(\tf)(x)= \int_{\gc'_\gt}\hf(\zeta)e^{-\zeta/x} d\zeta
\end{equation}
\noindent\begin{minipage}[c]{60mm}
\baselineskip=.8em

\hspace{.7cm}
\includegraphics[scale=1]{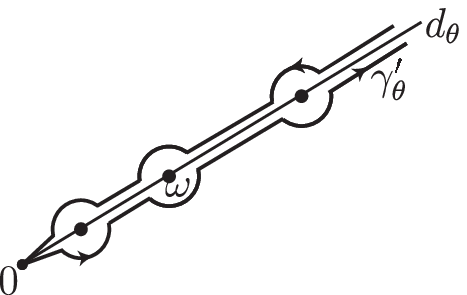}\\
\begin{footnotesize}
\baselineskip=.7em
  The two branches of $\gc'_\gt$ along $d_\gt$ are disjointly drawn to point out that they do not belong to the same sheet of the Riemann surface $\sR_\gO$.
  \end{footnotesize}
 \end{minipage}
\hskip 1cm
\begin{minipage}[c]{65mm}
\baselineskip=1.5em
where, as shown on the figure, $\gc'_\gt$ is a Hankel type path going along the straight line $d_\gt$ from infinity to 0 and back to infinity passing positively all singularities $\go\in\gO_\gt$ on both ways.\\ 

\end{minipage}\\  \\
The exponential growth of $\hf$ in direction $\theta$ guaranties the convergence of the integral for $x$ in a disc adherent to 0 with  $d_\gt$ as a diameter (a Borel disc) (\cf \cite{MR90} for example).

Without changing the value of the integral\footnote{ Contrarily to Formula (\ref{stomatbis}) which only requires the 1-summabilty of the series $\tf$, the individual resurgence and 1-summability  are not sufficient here. We do need  summable-resurgence.} the path $\gc'_\gt$ can be deformed 
\noindent\begin{minipage}[c]{60mm}
\vspace{2mm}

\hspace{1cm}
\includegraphics[scale=1]{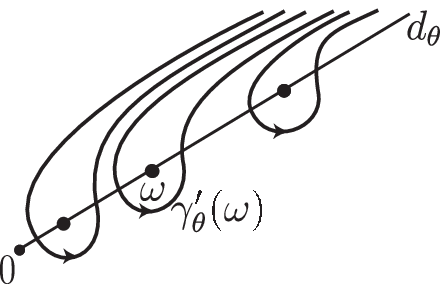
}
 \end{minipage}
\hskip 1cm
\begin{minipage}[c]{65mm}
\baselineskip=1.5em
into the union $\gc'_\gt=\cup_{\go\in\gO_\gt} \gc'_\gt(\go) $ of Hankel type paths $\gc'_\gt(\go)$  with asymptotic direction $\gt$ around each singular point $\go\in\gO_\gt$. \\

\end{minipage}\\ 

By means of translations from $\go$ to 0 and using the fact that holomorphic functions at $\go$ contributes 0 to the integral around $\go$ we can replace $\hf$ by its principal majors $\widecheck f^+(\go+\xi) $ at each $\go$ obtaining so


\noindent\begin{minipage}[c]{100mm}
\baselineskip=1.5em
\begin{equation}
s_{\gt^+}(\tf)-s_{\gt^-}(\tf)= \sum_{\go\in\gO_\gt} e^{-\go/x}\,\int_{\gc_\gt}\widecheck f ^+(\go+\xi)e^{-\xi/x} d\xi
\end{equation}
where, as shown on the figure, $\gc_\gt$ is a Hankel type path around  $0$ in direction $\gt$.
 \end{minipage}
\hskip 5mm
\begin{minipage}[c]{30mm}
\vspace{.5cm}
\qquad\includegraphics[scale=1]{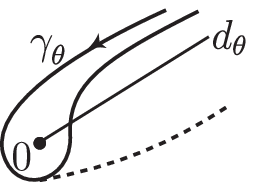}
\end{minipage}\\ \\


We claim that we can identify this linear combination of exponentials $\displaystyle e^{-\go/x}$ to the  linear combination in the right hand side of Formula (\ref{stomatbis}).

\begin{pro} Given $\go\in \gO_\gt$, the following identity holds for $\arg(x)\simeq \ugt$:
\begin{equation}\label{identity}
\int_{\gc_\gt}\widecheck f^+ (\go+\xi)e^{-\xi/x} d\xi=
s_{\gt^-}(\tF)(x)\, x^{L}\, \ddodp\, x^{-J_{1}}.
\end{equation}
\end{pro}
\proof The equality is obviously true when there is only one singular point $\go$ on the half line $d_\gt$. Assume then that there are at least two singular points on $d_\gt$.\\ The Borel transform $\hf(\xi)$ of $\tf(x)$ belongs to $\widehat{\sR es}_\gO^{\rm sum}$ and has simple moderate singularities. From Theorem \ref{thmsingf} we know that any major $\widecheck f(\go+\xi)$ at $\go$ is a polynomial in $\ln \xi$ with summable-resurgent coef\-ficients, possibly factored by complex powers $\xi^{\gl-1}$. The Laplace integral $\int_{\gc_\gt}\widecheck f^+ (\go+\xi)e^{-\xi/x} d\xi$ in the left hand side is then a polynomial in $\ln x$ whose coefficients are 1-Gevrey asymptotic functions on a germ of half plane bisected by the direction $\gt$. Deforming the path $\gc_\gt$ so as to move  the asymptotic direction to a neighboring direction $\gt+\gve$ we get 
1-Gevrey asymptotic coefficients on a sector of opening $ \rbrack \gt-\frac{\pi}{2} , \gt+\frac{\pi}{2}+\gve\lbrack$, that is, 1-sums of 1-summable series in direction $\displaystyle \gt+\frac{\gve}{2}$ (integer power series, possibly factored by non integer powers $x^\gl$).

The same property holds for the right hand side  and therefore, we can conclude by the  variant of  Watson Lemma  below. \qed

 \begin{lemma}[A variant of Watson Lemma]\label{Watson}\quad
 
Let $\gt\in\R/2\pi\Z$. 

Suppose that, to each $\go\in\gO_\gt$,  there is  a linear combination  $P_\go(\ln x)$ of polynomials in $\ln x$ with summable-resurgent coefficients in $\widetilde{\sR es}_{\gO-\go}^{\rm sum}$ possibly weighted by  complex powers of $x$, the powers being distinct modulo $\Z$.

Then, an identity $\displaystyle \sum_{\go\in\gO_\gt} P_\go(\ln x) \, e^{-\go/x} \equiv 0$ implies the nullity of each $P_\go$, \ie the nullity of each summable-resurgent coefficient in $P_\go$.
\end{lemma}

\proof By means of a rotation of the variable $x$ we can assume that $\gt=0$ so that all $\go$ are positive numbers. We range them by increasing order:
$ \go_1< \go_2<\cdots< \go_r$. If we factor $\displaystyle e^{-\go_1/x} $ the identity becomes 
$$
P_{\go_1}(\ln\, x) +\sum_{\ell=2}^r P_{\go_\ell}(\ln\, x)e^{-(\go_\ell-\go_1)/x} \equiv 0.
$$ 
Taking the asymptotic expansion at 0 on the positive real line of both sides of the identity proves that the asymptotic expansion of  $P_{\go_1}(\ln\, x)$ at 0 is 0. The same property holds for the coefficients of each power of $\ln\, x$. These coefficients are 1-summable powers series in which the sequence of exponents is a finite union of arithmetic sequences.  Watson Lemma applies to such a situation (\cf \cite[Thm 2.4.1.4.ii]{Mal95}) implying the unicity of the 1-sum. Hence, $P_{\go_1}\equiv 0$. The nullity of the other polynomials $P_\go$ follows by recursion. \qed

\begin{thm} {\bf Connection-to-Stokes  Formula}\\
Given $\gt\in\R/2\pi\Z$ an anti-Stokes direction,  the data of the Stokes-Ramis matrix  $C_{\ugt}$ and of the connection  matrix $K_{\ugt}^+$ are equivalent and the two matrices are related by the relation
\begin{center}
\begin{tabular}{| p{5.8cm}|}
\hline 
$\displaystyle \  \ C_{\ugt}=\int_{\gc_0} \frac{1}{\xi}\, \xi^L \,K_{\ugt}^+\, \xi^{-L}\,e^{-\xi}d\xi
$\tvi{20}{18}\\
\hline
\end{tabular}
\end{center}
where $\gc_0$ is  a Hankel type path around $\R^+$ run over by $\xi$ after the choice of the principal determination of its argument.
\end{thm}
\proof Note that $C_{\ugt}$ and $K_{\ugt}^+$ have the same block-structure with the same blocks of zeros and same arbitrary non-zero blocks. The map $K_{\ugt}^+ \mapsto C_{\ugt}$ is a linear map between two linear spaces of same dimension and to prove that it is bijective it is enough to prove that it is injective. 

Suppose $C_{\ugt}=0$; then, the two lateral sums $s_{\gt^+}(\tF)$ and  $s_{\gt^-}(\tF)$ glue together and the Borel transform $\hF(\xi)$ of $\tF(x)$ has no singularity on the line $d_\gt$. 
Hence, $K_{\ugt}^+=0$.

To prove the relation linking $C_{\ugt}$ to $K_{\ugt}^+$ it is sufficient to compare the first column-blocks $c_{\ugt}$ and $k_{\ugt}^+$ (\cf Section \ref{prepared}) and to prove that
\begin{equation*}
c_{\ugt}=\int_{\gc_0} \frac{1}{\xi}\, \xi^L\, k_{\ugt}^+\, \xi^{-J_1}\, e^{-\xi} d\xi
\end{equation*}
with the choice of the principal determination of argument along $\gc_0$.\\
 Equating the dominant terms in Identity (\ref{identity}) and using Theorem \ref{thmsingf} we obtain the new identity 
\begin{equation*}
\begin{array}{ccl}
x^L\, c_{\ugt}\, x^{-J_1}&=& \displaystyle \int_{\gc_\gt}u^L \,k_{\ugt}^+ \, u^{-J_1} e^{-u/x}du\\
\noalign{\medskip}
&=& \displaystyle  x^L\int_{\gc_0} \xi^L \,k_{\ugt}^+\, \xi^{-J_1} e^{-\xi}d\xi\, x^{-J_1} \quad(\text{setting } u=x\xi)\\
\end{array}
\end{equation*}
where for ${\rm arg}\, x= \ugt$, the path of integration $\gc_\gt$ has become the classical Hankel path $\gc_0$ around the non-negative real axis with argument from $-2\pi$ to 0. Hence the result. \qed

In restriction to the blocks attached to a given $\go$ we obtain the following statement.
\begin{cor}
For all $\go\in{\bf \gO}^*$, the data of $\dDodp$ and of $K_{(\ugo)}^+$ are equivalent and the two matrices are related by 
\begin{center}
\begin{tabular}{| p{5.8cm}|}
\hline 
$\displaystyle \  \ \dDodp=\int_{\gc_0} \frac{1}{\xi}\, \xi^L \,K_{(\ugo)}^+\, \xi^{-L}\,e^{-\xi}d\xi
$\tvi{20}{18}\\
\hline
\end{tabular}\, 
\end{center}
where $\gc_0$ stands for  a Hankel type path around $\R^+$ run over by $\xi$ after the choice of the principal determination of its argument.
\end{cor}

\bigskip

It can be useful, especially for effective numerical calculations to expand each entry of the Connection-to-Stokes Formula. The following corollary provides such an  expanded form. \\
Let $M^{(j,\ell);(k,r)}$  denote  the entry row $\ell $ of row-block $j$ and column $r$ in column-block $k$ of a $n\times n$-matrix $M$ split into blocks fitting the structure of  $L$. \\
\begin{cor}\label{CorConn2Stokes}
The various entries of the Connection-to-Stokes  Formula read 
\renewcommand{\arraystretch}{2} 
\begin{center}
\begin{tabular}{| p{8cm}|}
\hline 
$\displaystyle \  \ C_{\ugt}^{(j,\ell);(k,r)}= \sum_{p=0}^{n_j-\ell+r-1} 
\kappa_p(\gl_j-\gl_k) \, H_p^{(j,\ell);(k,r)}
$\tvi{25}{18}\\
\hline
\end{tabular}
\end{center}
where $\displaystyle \kappa_p(\gl_j-\gl_k)=2\pi i \frac{d^p}{dt^p}\Big(\frac{e^{-i\pi t}}{\gC(1-t)}\Big)\Big\vert_{ t=\gl_j-\gl_k}$ and \\
$\displaystyle H_p^{(j,\ell);(k,r)} = \sum_{\ell'-r'=p-r+1} (-1)^{r-1-r'} 
\frac{K_\gt^{+^{(j,\ell+\ell');(k,r'+1)}}}{(r-1-r')!\,\ell'!}\cdot\tvi{25}{0}$
\end{cor}

\subsection{Alien derivations}\label{aliender}

As already said in the introduction  of Section \ref{MeroInv}, the classifying set of meromorphic classes of formal gauge transformations of a given normal system  is naturally endowed with a structure of a Lie group. In this section, we explain how the dotted alien derivations can be defined as infinitesimal generators of this Lie group.
\bigskip

Alien derivations were first given a definition through the Borel plane by J. \'Ecalle (\cite{Ec81, Ec92}, \cf also \cite{Mal85,Sau05}). They have been defined as an average of analytic continuations in the Borel plane followed by a Laplace transform, the weights of the averaging being related to the Campbell-Hausdorf Formula. It is equivalent to see them as the homogeneous components of the logarithm of  graded Stokes automorphisms (\cf \cite{Ec81}, \cite[Th I.6.2]{Mal85}, \cite[Lemma 5 and pp. 35-38]{Sau05}). This is the viewpoint we adopt here, a viewpoint  already used in \cite{R89, MR90}. To perform it we need  
\begin{itemize}
\item to define an algebra where the Stokes automorphisms make sense and where the alien derivatives will live (to give them a chance of being ``derivatives'' we cannot keep working in vector spaces),
\item to define a graduation on the Stokes automorphisms.
\end{itemize}

Our aim being to define alien derivations of the solutions of System (\ref{A}) we proceed as follows.
\subsubsection{The algebra}
 
Consider the finitely many Stokes values $\gO=\{a_1,\dots,a_J\}$ associated with System (\ref{A}) and the $\Z$-module $\Z\gO$  they generate.  The lattice ${\Z\gO}$ may be dense in $\C$ (\cf the example of the hypergeometric equation $D_{13,1}$ below) and it contains  ${\bf \gO}=\{a_j-a_\ell\}$. 
Simultaneously, consider the set $\gL=\{\gl_1,\dots,\gl_J\}$ of exponents of formal monodromy of System (\ref{A}) and the lattice $\Z\gL$ it generates.
\smallskip

We introduce the set $\stH_{\Z\gL,\Z\gO}$ of all series solutions of linear differential equations with meromorphic coefficients, single level 0 or 1, exponents of formal monodromy in $\Z\gL$ and Stokes values in $\Z\gO$. Note that all (convergent) meromorphic series at 0 belong to
$\stH_{\Z\gL,\Z\gO}$.

\begin{pro}
The set $\stH=\stH_{\Z\gL,\Z\gO}$ is a differential sub-algebra of the algebra of all formal meromorphic series at $x=0$.\\
Consequently, its counter-part $\widehat\sH$ in the Borel plane is a convolution algebra where multiplication by $1/\xi$ is the derivation; its elements have simple-moderate singularities at finitely many points of $\Z\gO$.
\end{pro}
\proof The set $\stH$ is stable with respect to sums and products. Indeed, given two series solutions of linear differential equations $D_1y=0$ and $D_2y =0$ respectively, their sum satisfies a linear equation $Dy=0$ where $D$ is the left smaller common multiple of $D_1$ and $D_2$, the exponents of formal monodromy and the Stokes values being the union of those of the initial equations. Their product satisfies a linear differential equation $\gD y=0$ where $\gD$ is the symmetric tensor product of $D_1$ and $D_2$ and where the exponents of formal monodromy and the Stokes values are the two by two sums  of those of the initial equations.\\
The set $\stH$ is stable under derivation with respect to $x$ since the derivative of an element $\tf$ of $\stH$ satisfies an equation of the same type obtained by conveniently derivating a differential equation satisfied by $\tf$.
\qed

To prove directly in the Borel plane that $\widehat\sH$ is a convolution algebra is not so easy. Given $\hf$ and $\hg$ in $\widehat\sH$ with finite singular supports $\gO(\hf)$ and $\gO(\hg)$ respectively,  the main point is to prove that the convolution $\hf*\hg$ which is well defined near  0 can be continued to the whole Riemann surface $\sR_{\gO(\hf) +\gO(\hg)}$. To this end, one could generalize the technique of two 
intertwined combs to build $\sR$-symmetrically contractile paths (\cf \cite[Lemma 3, Figures 5 and 10]{Sau05}): one should start here with a comb with upwards nails at ${\bf\gO}(\hf)$, a comb with downwards nails at  $-{\bf\gO}(\hg)$, symmetric of  ${\bf\gO}(\hg)$ with respect to 0, and an elastic rope tied to the two nails 0 to materialize the convenient paths of analytic continuation when one moves the second comb. We won't formalize such a proof. Exponential growth at infinity and simple-moderate singularities are preserved by convolution.
\medskip

In order to define alien derivatives of the entries of the gauge transformation $\tF(x)$ we limit ourselves to consider the differential sub-algebra 
$\stH_{\tF}$ of $\stH$ generated by the entries of $\tF(x)$.
Formula (\ref{stoformel}) as noticed at the end of Section \ref{formalsto} shows that the Stokes automorphism $St_\gt$  does not act inside $\stH_{\tF}$ since the image of a formal series may involve complex powers of $x$, logarithms and exponentials. 
It is thus natural to extend the differential algebra $\stH_{\tF}$ into the polynomial algebra of ``resurgent symbols''\footnote{For non-linear situations one has to consider resurgent symbols with infinitely many exponential terms.}
$$
\stH_{\tF \tY_0}=\stH_{\tF}\,\big\lbrack(x^{\gl})_{\gl\in\Z\gL},\ln x, (e^{-\go/x})_{\go\in{\Z\gO}}\big\rbrack.
$$
This is the differential algebra we are willing to work in.  The coefficients are formal 1-summable series, $x^{\gl},\,\ln x, \,e^{-\go/x}$ are formal indeterminates  satisfying the usual rules. The derivation is $d/dx$.
\subsubsection {Extended Stokes automorphisms}  \label{extauto}

Let $\gt\in\R/2\pi\Z$ be an anti-Stokes direction associated with System (\ref{A}) and $\ugt\in\R$ its principal determination.\\

In Section \ref{formalsto},  we described the Stokes automorphism as a map acting on formal solutions by means of a choice of a determination of the argument $\ugt$. Such a definition can be extended  into an automorphism of the differential algebra $\stH_{\tF \tY_0}$ as follows.
\goodbreak

\begin{pro}[Extended Stokes automorphism]\quad\\
The  Stokes automorphism $St_{\ugt}$ can be extended into a differential  unipotent automorphism, still denoted $St_{\ugt}$, of the differential algebra $\stH_{\tF \tY_0}$  by setting:
\begin{description}
\item{-- }  power factors $x^{\gl }$,  logarithms and exponentials are kept fixed,
\item{-- }   $\tF$ is changed into $\displaystyle \tF\, x^L\,e^{Q(1/x)}(I_n+C_{\ugt})\,e^{-Q(1/x)}\,x^{-L}\tvi{15}{0}$. 
\end{description}
\end{pro}
\proof 
Any resurgent symbol of $\stH_{\tF \tY_0}$ has a  unique expression of the form 
$\displaystyle \sum \tg_j\,x^{\gl_j}\,\ln^{p_j} x\, e^{-\go_j/x} $
the sum running on finitely many distinct triples $(\gl_j,p_j,\go_j)\in \Z\gL\times \N\times \Z\gO$ if we assume, in addition, that $0\leq \Re \gl_j<1$ for all $j$ (and $\tg_j\neq 0$). It can be isomorphically sent to the actual resurgent symbol    
$\displaystyle \sum s_{\gt^-}(\tg_j)(x)\,x^{\gl_j}\,\ln^{p_j}(x)\, e^{-\go_j/x} $ for $\displaystyle \arg x\in \big\rbrack \ugt-\frac{\pi}{2}, \ugt+\frac{\pi}{2} \big\lbrack$. With this isomorphism the map $St_{\ugt}$ as defined in the proposition reads  as the map 
\begin{equation*}
\sum s_{\gt^-}(\tg_j)(x)\,x^{\gl_j}\,\ln^{p_j}(x)\, e^{-\go_j/x} \mapsto
\sum s_{\gt^+}(\tg_j)(x)\,x^{\gl_j}\,\ln^{p_j}(x)\, e^{-\go_j/x}
\end{equation*}
in the space of actual resurgent symbols of the form 
$$\displaystyle \sum s_{\gt^-}(\tg_j)(x)\,x^{\gl_j}\ln^{p_j}(x) e^{-\go_j/x} 
\text{ with } \tg_j\in\stH_{\tF \tY_0} \text{ and }  \arg x\in \big\rbrack \ugt-\frac{\pi}{2}, \ugt+\frac{\pi}{2} \big\lbrack$$ 
if one expands $s_{\gt^+}(\tg_j)(x)$ according to Formula (\ref{stomatbis}) (Since  $\tg_j$ is 1-summable it has a well-defined Stokes phenomenon and the result does not depend on the way $\tg_j$ is expanded in terms of the entries of $\tF$). 
Since summations $s_{\gt^-}$ and $s_{\gt^+}$ are automorphisms of differential algebras so is $St_{\ugt}$. The reciprocal map $St_{\ugt}^{-1}$  is obtained by keeping the same space of actual resurgent symbols while exchanging the roles of $s_{\gt^-}$ and $s_{\gt^+}$ and then the Stokes matrix $I_n+C_{\ugt}$ by its inverse $(I_n+C_{\ugt})^{-1}$.

\bigskip

 The extended Stokes automorphism keeps being unipotent. Indeed, this results from the fact that it is already unipotent when acting in the space of formal solutions of System (\ref{A}) (its matrix $I_n+C_{\ugt}$ is unipotent). This can also be seen as follows:
 if $\tF_{j,k}$ is an entry of $\tF$ then 
$St_\gt(\tF_{j,k})(x)=\tF_{j,k}(x)+\sum g_{\gl,\ell,\go}(x) \, x^\gl \ln^{\ell}(x)\, e^{-\go/x}$
 where the sum runs on finitely many $\gl$ and $\ell$ and  finitely many $\go$ in ${\bf\gO}_\gt$. The coefficients $g_{\gl.\ell.\go}$ are themselves elements of $\stH_{\tF \tY_0}$. The Stokes values of $g_{\gl.\ell.\go}$ are among those of $\tF$ translated by $-\go$ hence none is left on the half-line $d_\gt$ after finitely many applications of $St_\gt$. \qed

With this extended definition we can  now write \begin{equation}
St_{\ugt}(\tF(x)) =\tF(x)\, x^L\, e^{Q(1/x)} (I_n+C_{\ugt})e^{-Q(1/x)}\, x^{-L}
\end{equation}
 instead of $\displaystyle  
St_{\ugt}(\tF(x)\, x^L\, e^{Q(1/x)}) = \tF(x)\, x^L\, e^{Q(1/x)} (I_n+C_{\ugt})$ only (\cf  (\ref{stoformel})).  

 \subsubsection{Graduation on the Stokes automorphisms}
 
 The graduation is built so as to discriminate between the different sub-matrices $\dD_{\ugo}^+$ of the Stokes matrix $C_{\ugt}$ (\cf Section \ref{MapDelta+}). This is done using the  {\it exponential torus}  $\sT$ of System (\ref{A}).\\

 Here is how to define  $\sT$ (\cf \cite{MR90,MR91}):
  Let $b_1,\dots,b_\nu$ be a basis of the lattice ${\Z \gO}$.
 The polynomials $p_1(1/x)=-b_1/x,\dots,p_\nu(1/x)=-b_\nu/x$ form a basis of the  $\Z$-module generated by the determining polynomials $q_1, q_2,\dots,q_n$ of the diagonal of $Q(1/x)$. The exponential torus takes into account the indetermination  of a formal exponential $\displaystyle e^p$ by associating with $\displaystyle e^p$ its complex multiples $\gl \displaystyle e^p$ (see Footnote \ref{exptorus}). To define the exponential torus $\sT$ one introduces $\nu$ indeterminates $\underline{\gl}=(\gl_1,\gl_2,\dots,\gl_\nu)$ and associates $\displaystyle \gl_j e^{p_j}$ with each $\displaystyle  e^{p_j}$. One can extend $\sT$ to the algebra $\tH_{\tF \tY_0}$ by letting it act trivially on power series, complex powers of $x$ and logarithms. We obtain thus a family of automorphisms of $\stH_{\tF \tY_0}$ with $\nu$ parameters.\\
  For $j=1,2,\dots,n$, denote $\underline{m}_j=(m_{j,1}, m_{j,2},\dots,m_{j,\nu})$ the components of the $q_j$'s with respect to the $p_r$'s.   With respect to the formal fundamental solution $\displaystyle \tF(x)\,x^L\,e^{Q(1/x)}$ the exponential torus  has a matrix representation of the  form
\begin{equation}
T_{\underline{\gl}}=\diag(\underline{\gl}^{\underline{m}_1},\underline{\gl}^{\underline{m}_2},\dots,\underline{\gl}^{\underline{m}_n})
\end{equation}
where the notation $\underline{\gl}^{\underline{m}}$ stands for the product 
$\gl_1^{m_1}\,\gl_2^{m_2}\dots\gl_\nu^{m_\nu}$. Its action on the Stokes automorphism $St_\ugt$  generates a group $G_{\underline\gl}$ of unipotent   matrices $I_n+T_{\underline{\gl}}\, C_{\ugt}\,T_{\underline{\gl}}^{-1}$ with $\nu$ parameters. 
These  matrices are  polynomials in the parameters $\gl_r$ and their inverses. Keep denoting ${\bf \gO}_\gt$ the set of non-zero Stokes values $a_j-a_k$ on the line $d_\gt$ (\cf Section \ref{knownresults}).
To each Sto\-kes value $\go\in{\bf \gO}_\gt$ there is a unique collection $\underline{m}(\go)=(m_1(\go),\dots,m_\nu(\go))$ of weights such that $\displaystyle \sum_{r=1}^\nu m_r(\go)\, p_r(x)=-\frac{\go}{x}$. One can check that the matrix  $\dD_{\ugo}^+$ (\cf Section \ref{MapDelta+}) is the coefficient of the monomial   $\underline{\gl}^{\underline{m}(\go)}$ in $T_{\underline{\gl}}\, C_{\ugt}\,T_{\underline{\gl}}^{-1}$:
\begin{equation}
T_{\underline{\gl}}\, C_{\ugt}\,T_{\underline{\gl}}^{-1} = \sum_{\go\in{\bf\gO}_\gt} \dD_{\ugo}^+ \, \underline{\gl}^{\underline{m}(\go)}.
\end{equation}
\subsubsection{Definition of the alien derivations}

As a unipotent graded group the group $G_{\underline\gl}$ admits  infinitesimal generators $\dD_{\ugo}$ in the sense that 
\begin{equation}
I_n+T_{\underline{\gl}}\, C_{\ugt}\,T_{\underline{\gl}}^{-1} = 
\exp\Big(\sum_{\go\in{\bf\gO}_\gt} \dD_{\ugo}\, \underline{\gl}^{\underline{m}(\go)}\Big).
\end{equation}


\begin{dfn}[dotted and undotted alien derivations]\quad 

{\rm \begin{description}
\item{$\bullet$} The {\it dotted alien derivations} $\dD_{\ugo}$ are the transformations with matrix the coefficient of  $\underline{\gl}^{\underline{m}(\go)}$ in the expansion of the logarithm $\ln \big(I_n+T_{\underline{\gl}}\, C_{\ugt}\,T_{\underline{\gl}}^{-1}\big)$ (matrix in the chosen basis of formal solutions $\displaystyle \tF(x)\,x^L\,e^{Q(1/x)}$).
\begin{equation}
\ln \big(I_n+T_{\underline{\gl}}\, C_{\ugt}\,T_{\underline{\gl}}^{-1}\big)=\sum_{\go\in{\bf\gO}_\gt} \dD_{\ugo}\, \underline{\gl}^{\underline{m}(\go)}
\end{equation}
\item{$\bullet$} The {\it alien derivation} $\gD_{\ugo}$ is defined by $\gD_{\ugo}= e^{+\go/x}\dD_{\ugo}$.

\end{description}
}\end{dfn}

When $\go$ does not belong to ${\bf\gO}_\gt$ for any $\gt$ the alien derivations $\dD_{\ugo}$ and $\gD_{\ugo}$ are equal to 0. 
\\
Like $St_{\ugt}^+$  the dotted alien derivations act trivially on $x^L$ and $\displaystyle e^{Q(1/x)}$ so that 
$$\displaystyle \dD_{\ugo}\big(\tF(x)\big)=\tF(x)\,x^L e^{-\go/x}\dD_{\ugo}\,x^{-L}.$$ 
and 
$$\displaystyle \gD_{\ugo}\big(\tF(x)\big)=\tF(x)\,x^L \dD_{\ugo}\,x^{-L}.$$ 
The alien derivations are derivations by construction and they   commute with  the usual derivation $d/dx$.
\bigskip

We end this section with a remark on the various choices made.\\
We saw that the meromorphic classifying set is given by the Stokes automorphisms $St_{\gt,\tF}$ defined, for all anti-Stokes direction, as gauge transformations of the normal form. To look at the Stokes automorphisms $St_{\ugt}$ defined as linear maps on the space of formal solutions of System (\ref{A}) and get their extended forms we needed to choose a determination of the argument and an actual form of the formal exponentials. With different choices the Stokes matrix is conjugate under the iterated action of the formal monodromy with matrix $\widehat{M}=x^{2\pi i L}$ with respect to $\tY_0=x^L\,e^{Q(1/x)}$
 and under the action of the exponential torus $\sT$. The same conjugacy actions must be taken into account when performing the analytic classification with alien derivations.
\subsubsection{Bridge equation}

 The ``definition'' formula above  rewritten  in the form 
\begin{equation}
 \gD_{\ugo}\big(\tF(x)x^L\big)= \tF(x)\,x^L\dD_{\ugo}
\end{equation} 
can be seen as \'Ecalle's Bridge Equation. The name ``bridge'' comes from the fact that  the equation links alien derivatives (left hand side of the bridge equation) to ordinary derivatives (right hand side of the bridge equation). Indeed, the right hand side can be seen as an ordinary derivative as follows.
We consider only the first $n_1$ columns, the calculation for the $k^{th}$ block of columns  being the same after multiplication by $\displaystyle e^{a_k/x}$. Introduce the general solution (also said formal integral) of System (\ref{A}) which has the form 
$$
\tf(x)\,x^{L_1}M_0 +\sum_{\go\in\{a_2,\dots,a_J\}} \phi_\ugo(x) \,M_\ugo \,e^{-\go/x}
$$
where the $\phi$'s are formal-log series and the $M$'s are arbitrary constant matrices. The alien derivative of $\displaystyle \tf(x)\,x^{L_1}$ at $\go$  is the derivative of this general solution with respect to $\displaystyle e^{-\go/x}$ considered as an independent variable\footnote{ In the case of scalar solutions of an equation instead of a system, J. \'Ecalle takes derivatives with respect to the constants. Since, here, the constant  coefficients of the various exponentials are matrices it is more convenient to derivate with respect to the exponentials themselves.}, with a convenient choice of the matrices $M$ (\'Ecalle's analytic invariants $\dD_\ugo$).

 In \'Ecalle's approach, alien derivatives are defined as an average of analytic continuations in the Borel plane followed by a Laplace transform. The Bridge Equation results from the fact that dotted alien derivatives commute with the derivation $d/dx$. In this approach, it is nothing more than the definition formula.

\protect\boldmath
\subsection{The example of the generalized hypergeometric equation  $D_{13,1}$}\protect\unboldmath

We consider the generalized hypergeometric equation of order 13
\begin{equation}
{\bf D_{13,1}}({\bf y})\equiv \Big(x\frac{d{\bf y}}{dx}-\mu {\bf y}\Big) - x \prod_{j=1}^{13}\Big(x\frac{d}{dx}-(\nu_j-1)\Big){\bf y}=0
\end{equation}
where $\mu$ and the $\nu_j$'s are complex parameters\footnote{The irregular singular point, usually put at infinity (\cf \cite{DM89}),  is here located at 0.}. Its Newton polygon at 0 has a slope 0 of length 1 and a slope 1/12 of length 12. Putting $x=t^{12}$ the equation becomes 
\begin{equation}\label{hyperg}
 D_{13,1}(y)\equiv \Big( \frac{1}{12}\,t\,\frac{dy}{dt}-\mu y\Big)-t^{12}\prod_{j=1}^{13}\Big(\frac{1}{12}\, t\, \frac{d}{dt}-(\nu_j-1)\Big)\,y=0.
\end{equation}
We keep using the equation itself taking benefit of  having a quite simple equation but we could as well commute to the companion system.

The hypergeometric equation (\ref{hyperg}) is of single level 1. Its determining polynomials $q_1, q_2,\dots,q_{13}$ are calculated in \cite{DM89} and are
\renewcommand{\arraystretch}{1.5} 
$$\begin{array}{llll}
q_1=0 &\  q_2\ =-12/t&\quad q_3\ =-12\gz/ t&\quad q_4\ =-12\gz^2/ t\\ &\ q_5\ =-12\gz^3/ t &\quad q_6\ =-12\gz^4/ t&\quad q_7\ = -12\gz^5/ t\\
&\  q_8\ =-12\gz^6/ t&\quad q_9\ =-12\gz^7/ t &\quad  q_{10}=-12\gz^8/ t\\
&\ q_{11}=-12\gz^9 /t&\quad q_{12}=-12\gz^{10}/ t&\quad q_{13}=-12\gz^{11} /t\\
\end{array}$$
 where $\gz$ stands for the twelfth primitive root of unit $\displaystyle \gz=e^{2\pi i/12}$. A formal fundamental solution $\tF(t)\, t^{L}e^{Q(1/t)}$ reads 
 \begin{equation}
 \displaystyle\begin{bmatrix}
\tF^1(t)\,t^{12\mu}\, e^{q_1(1/t)}&\tF^{2}(t)\,t^{-
12\gl}\, e^{q_{2}(1/t)} &\cdots&\tF^{13}(t)\,t^{-12
\gl}\, e^{q_{13}(1/t)}
\end{bmatrix}
\end{equation}
 where $ \displaystyle \gl=\frac{1}{12}\Big(\frac{13}{2}+\mu-\sum_{j=1}^{13}\nu_j  \Big)$.
 
 \noindent\begin{minipage}[c]{57mm}
\begin{center}
\vspace{-.3cm}
\includegraphics[scale=1]{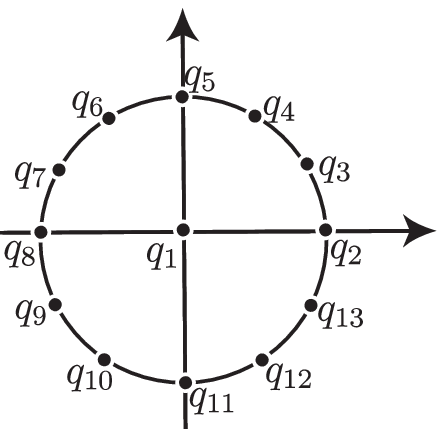}
\end{center}
 \end{minipage}
\hfill\begin{minipage}[c]{78mm}
\baselineskip=1.5em
\bigskip

From the fact that the minimal polynomial of $\gz$ is the cyclotomic polynomial $X^4-X^2+1$ we deduce that a basis for the lattice built on the $q_j$'s is given by 
$$
\begin{array}{ll}
p_1=q_2=-12/ t&\quad p_2=q_3=-12\gz/ t\\
p_3=q_4=-12\gz^2/ t&\quad p_4=q_5\ =-12\gz^3 /t\\
\end{array}
$$
\end{minipage}\\ \\
Note that the lattice built on the coefficients $12, 12\gz,12\gz^2,12\gz^3$ is dense in $\C$. However, only finitely many values are Stokes values for $D_{13,1}$. \\
In the $p$-basis the determining polynomials read 
\renewcommand{\arraystretch}{1.5} 
$$\begin{array}{lllll}
q_1=0 &\  q_2\ =p_1&\quad q_3\ =p_2&\quad q_4\ =p_3&\quad q_5\ =p_4\\
 &\  q_6\ =p_3-p_1&\quad q_7\ =p_4-p_2&&\\
 &\  q_{8}\ =-p_1&\quad q_{9}\ =-p_2&\quad q_{10}=-p_3&\quad q_{11}=-p_4\\
 &\ q_{12}=-p_3+p_1&\quad q_{13}=-p_4+p_2&&\\
\end{array}$$
so that, denoting ${\underline \gl}=(\gl_1,\gl_2,\gl_3,\gl_4)$, the matrix of the exponential torus is given by 
$$
T_{\underline\gl}=\diag\Big(1,\gl_1,\gl_2,\gl_3,\gl_4,\frac{\gl_3}{\gl_1}, \frac{\gl_4}{\gl_2},\frac{1}{\gl_1},\frac{1}{\gl_2},\frac{1}{\gl_3},\frac{1}{\gl_4},\frac{\gl_1}{\gl_3},\frac{\gl_2}{\gl_4} \Big).
$$
Let us consider the anti-Stokes direction $\gt=0$ with principal determination $\ugt=0$ and denote $C_0$ the corresponding Stokes matrix. \\
Let $E_{j,\ell}$ denote the elementary $13\times 13$-matrix the entries of which are all 0 but the one at row $j$ and column $\ell$ which is equal to 1.\\
The Stokes values $\go$ supported by the half line $d_0$, issued from 0 with argument 0 (\ie $\go\in{\bf \gO}_0$), are  
$$\renewcommand{\arraystretch}{1.5} 
\begin{array}{l}
\go=12 \text{ associated with } q_1-q_8=q_2-q_1=q_4-q_6=q_{12}-q_{10}=p_1\\
 \go=12 \sqrt{3}\text{ associated with } q_3-q_7=q_{13}-q_9=\sqrt{3}\,p_1 \\
 \go=24 \text{ associated with } q_2-q_8=2\,p_1\\
\end{array}$$
and consequently, the Stokes matrix satisfies 
\renewcommand{\arraystretch}{1.8} 
$$
\begin{array}{ll}
T_{\underline \gl} C_0 T_{\underline \gl}^{-1}=&\ \ \Big(c^{(1,8)}E_{1,8}+c^{(2,1)}E_{2,1}+c^{(4,6)}E_{4,6}+c^{(12,10)}E_{12,10}\Big) \,\gl_1\\
 &+ \Big(c^{(3,7)}E_{3,7} +c^{(13,9)}E_{13,9}  \Big)\, \displaystyle \frac{\gl_2^2}{\gl_4}\\
 &+ \Big(c^{(2,8)}E_{2,8}  \Big)\,\gl_1^2\\
\end{array}$$
since one has also the relations $q_3-q_7=q_{13}-q_9= 2p_2-p_4$.
In other words, we obtain:
\renewcommand{\arraystretch}{1.8} 
$$
\begin{array}{lcl}
\dD_{12}^+ &=& c^{(1,8)}E_{1,8}+c^{(2,1)}E_{2,1}+c^{(4,6)}E_{4,6}+c^{(12,10)}E_{12,10},\\
 \dD_{12\sqrt{3}}^+&=&c^{(3,7)}E_{3,7} +c^{(13,9)}E_{13,9},\\
 \dD_{24}^+&=&c^{(2,8)}E_{2,8}.\\
\end{array}$$
Note that $T_{\underline \gl} C_0 T_{\underline \gl}^{-1}$ does not depend on the parameter $\gl_3$ which, in turn, would appear in $T_{\underline \gl} C_{-\pi} T_{\underline \gl}^{-1}$.
The Stokes multipliers $c^{(j,\ell)}$ are made explicit in term of Barnes integrals in \cite{DM89}.\\
We know from the previous section that  the alien derivations in the various Stokes values belonging to ${\bf \gO}_0$ are given in terms of the Stokes multipliers above by taking the logarithm of $I_{13}+T_{\underline \gl} C_0 T_{\underline \gl}^{-1}$. From the relations \\
$(\dD_{12}^+)^2=c^{(2,1)}\,c^{(1,8)}\,E_{2,8}$ and\\ $\dD_{12}^+ \dD_{12\sqrt{3}}^+= \dD_{12\sqrt{3}}^+ \dD_{12}^+= \dD_{12}^+\dD_{24}^+=\cdots = (\dD_{24}^+)^2=0\tvi{20}{0}$ we obtain 

\begin{equation*}
\renewcommand{\arraystretch}{2} 
\begin{array}{l}
\ln\big(I_{13}+T_{\underline \gl} C_0 T_{\underline \gl}^{-1}\big)\\
\displaystyle\qquad =
\Big( \dD_{12}^+\gl_1+\dD_{12\sqrt{3}}^+\frac{\gl_2^2}{\gl_4}+\dD_{24}^+\gl_1^2\Big)-\frac{1}{2} \Big( \dD_{12}^+\gl_1+\cdots\Big)^2 +\cdots\\
\displaystyle\qquad = \ \ \dD_{12}^+\gl_1+\dD_{12\sqrt{3}}^+\frac{\gl_2^2}{\gl_4}+\Big(\dD_{24}^+-\frac{1}{2}(\dD_{12}^+)^2\Big)\gl_1^2.\\
\end{array}
\end{equation*}
Hence,
\begin{equation}\label{average}
\begin{tabular}{| p{12cm}|}
\hline 
$\displaystyle
\renewcommand{\arraystretch}{2} 
\begin{array}{lcl}
 \dD_{12}&=&\dD_{12}^+\ =\ c^{(1,8)}E_{1,8}+c^{(2,1)}E_{2,1}+c^{(4,6)}E_{4,6}+c^{(12,10)}E_{12,10}\\
\dD_{12\sqrt{3}}&=&\dD_{12\sqrt{3}}^+\ \ =\ c^{(3,7)}E_{3,7} +c^{(13,9)}E_{13,9}\\
\dD_{24}&=&\displaystyle \dD_{24}^+-\frac{1}{2}(\dD_{12}^+)^2 \ = \ 
\Big(c^{(2,8)}-\frac{1}{2} c^{(2,1)}\,c^{(1,8)} \Big)\,E_{2,8}\tvi{0}{20}\\
\end{array}
$\\ 
\hline
\end{tabular}
\end{equation}
 
 We can develop these formul\ae\ by writing $\displaystyle \dD_{12}(\tF)=\tF\, t^L\,e^{-12/t}\,\dD_{12}\, t^{-L}$ or $\displaystyle \gD_{12}(\tF)=\tF\,t^L\, \dD_{12}\, t^{-L}$  and so on\dots, \ie
\begin{equation}\label{HyperBridge}
\left\lbrace
\begin{array}{ll}
\gD_{12}(\tF^1)=\tF^2\,t^{-12(\gl+\mu)}\,c^{(2,1)}\qquad&\gD_{12}(\tF^6)=\tF^4\,c^{(4,6)}\\
\gD_{12}(\tF^8)=\tF^1\,t^{12(\gl+\mu)}\,c^{(1,8)}\qquad &\gD_{12}(\tF^{10})=\tF^{12}\,c^{(12,10)}\\
\gD_{12\sqrt{3}}(\tF^7)=\tF^3\,c^{(3,7)}\qquad &\gD_{12\sqrt{3}}(\tF^9)=\tF^{13}\,c^{(13,9)}\\
\displaystyle \gD_{24}(\tF^8)=\tF^2\,\Big(c^{(2,8)}-\frac{1}{2}c^{(2,1)}c^{(1,8)}\Big)&\\
\end{array}\right.
\end{equation}
all other alien derivatives on the real positive line being 0.
 
 Equations (\ref{HyperBridge}) can be seen as {\it \'Ecalle bridge equations}.
 \bigskip
 
 \begin{rem} {\rm We end this example with a comment on, for instance, the last formula in (\ref{average}) compared to those of \cite[Lemma 5]{Sau05}  deduced directly from the definition of alien derivations by analytic continuation in the Borel plane. If we consider that there is on $d_0=\R^+$ the three  singularities $\go_1=12,\go_2=12\sqrt{3}$ and $\go_3=24$ then our formula does not fit the expression given in Lemma 5 of \cite{Sau05} for $\gD_{\go_3}$.  To fit it in, we should actually re-introduce the  missing singular point $\go_1=24-12\sqrt{3}$ attached to the monomial $\displaystyle \frac{\gl_1^2\gl_4}{\gl_2^2}$ ---although with a null coefficient $\dD_{\go_1}^+=0$--- which combined with the monomial $\displaystyle \frac{\gl_2^2}{\gl_4}$ attached to $12\sqrt{3}$ gives the monomial $\gl_1^2$ attached to 24.
  Set $\go_2=12, \go_3=12\sqrt{3}, \go_4=24$ and the calculation, in both cases, gives the same result 
\begin{center}
$ \begin{array}{l}
\displaystyle \gD_{\go_4}= \gD_{\go_4}^+ -\frac{1}{2}({\gD_{\go_2}^+}^2 + \gD_{\go_1}^+\gD_{\go_3}^+ +\gD_{\go_3}^+ \gD_{\go_1}^+)\\
\displaystyle \hspace{2cm} +\frac{1}{3} ( {\gD_{\go_1}^+}^2\gD_{\go_2}^+ +\gD_{\go_1}^+\gD_{\go_2}^+\gD_{\go_1}^+ +\gD_{\go_2}^+{\gD_{\go_1}^+}^2) -\frac{1}{4}{\gD_{\go_1}^+}^4\\
\end{array}$
 \end{center} 
if one takes into account the fact that $\gD_{\go_1}^+=0$.}
 \end{rem}
  
 \subsection{An example with resonance}
 We consider the  system
\renewcommand{\arraystretch}{1} 
 \begin{equation*}
x^2\frac{dY}{dx}= \begin{bmatrix}
0&0&0&0\\
x^2&1&x&0\\
x^2&0&1&x\\
x^2&0&0&1\\
\end{bmatrix}\,Y \leqno{(\mathfrak S)}
\end{equation*}
 and its formal fundamental solution $\tY(x)=\tF(x)\,x^L\,e^{Q(1/x)}$ 
 where 
 \begin{itemize}
\item $Q(1/x)=\diag(0,-1/x,-1/x,-1/x)$, (Hence, the system has the unique level 1 and the Stokes values $\pm 1$),
\item $\displaystyle L=\begin{bmatrix}
0&0&0&0\\
0&0&1&0\\
0&0&0&1\\
0&0&0&0\\
\end{bmatrix}$ ($L$ is not diagonal; hence the resonance),
\item $\displaystyle \tF(x)=\begin{bmatrix}
1&0&0&0\\
\tf_2&1&0&0\\
\tf_3&0&1&0\\
\tf_4&0&0&1
\end{bmatrix}$ 
is a power series satisfying $\tF(x)=I_4+O(x^2)$.
\end{itemize}
The system admits the two anti-Stokes directions $\gt=0$ and $\gt=\pi$. Obviously, the Stokes matrix in direction $\pi$  is trivial: $I_4+C_{\pi}=I_4$.

We consider the anti-Stokes direction $\gt=0$ supporting the unique  Stokes value $\go=1$. Our aim is  the calculation of the alien derivation $\gD_1$  in terms of the Stokes multipliers in direction $\gt=0$.  Actually, although System $(\mathfrak S)$ is quite a little bit involved since it exhibits resonance, it is simple enough to allow an exact calculation of the Stokes multipliers as below. We will then be able to give an exact calculation for $\gD_1$.

One can check that  the series $\tf_j$'s are the unique solutions of the system \renewcommand{\arraystretch}{2.2} 
\begin{equation}
\left\lbrace
\begin{array}{ccl}
\displaystyle x^2\frac{d\tf_2}{dx}-\tf_2&=&x^2+x\tf_3\\
\displaystyle x^2\frac{d\tf_3}{dx}-\tf_3&=&x^2+x\tf_4\\
\displaystyle x^2\frac{d\tf_4}{dx}-\tf_4&=&x^2\\
\end{array}\right.
\end{equation}
satisfying the condition $\tf_j(x)=O(x^2)$.
\renewcommand{\arraystretch}{2} 
It results that their Borel transforms $\hf_j$ are given by
\begin{equation}
\left\lbrace
\begin{array}{ccl}
\displaystyle \hf_2(1+\xi)&=&\displaystyle \frac{1}{\xi}\Big(\frac{1}{2}(6-\pi^2+4\pi i)+ (2+\pi i)\ln \xi +\frac{1}{2}\ln^2 \xi  \Big) +3\\
\displaystyle  \hf_3(1+\xi)&=&\displaystyle \frac{1}{\xi}\Big((2+\pi i)+\ln\xi\Big)+2\\
\displaystyle  \hf_4(1+\xi)&=&\displaystyle \frac{1}{\xi}+1\\
\end{array}\right.
\end{equation}
and consequently, the connection matrix $K_{(1)}^+$ is given by
\begin{equation*}
K_{(1)}^+= \begin{bmatrix}
0&0&0&0\\
k_2=\frac{1}{2}(6-\pi^2+4\pi i)&\ 0\ &\ 0\ &\ 0\ \\
k_3=2+\pi i&0&0&0\\
k_4=1&0&0&0\\
\end{bmatrix}\, .
\end{equation*}
From Corollary \ref{CorConn2Stokes} we deduce that the Stokes multipliers $C_0^{(2,1)}, C_0^{(3,1)}$ and $C_0^{(4,1)}$ are 
\begin{equation*}
\left\lbrace
\begin{array}{ccl}
C_0^{(2,1)}&=& \gk_0(0)\,k_2 +\gk_1(0)\,k_3+\frac{1}{2} \gk_2(0)\,k_4\\
C_0^{(3,1)}&=& \gk_0(0)\,k_3+\gk_1(0)\,k_4\\
C_0^{(4,1)}&=& \gk_0(0)\,k_4\\
\end{array}\right.
\end{equation*}
Recall that $\displaystyle \gk_p(\gl)=2\pi i\frac{d^p}{dt^p}\Big(\frac{e^{-i\pi t}}{\gC(1-t)}\Big)_{\Big\vert_{t=\gl}}$ and then,
\begin{equation*}
\left\lbrace
\begin{array}{ccl}
\gk_0(0)&=& 2\pi i\\
\gk_1(0)&=& 2\pi^2-2\pi i\gc\\
\gk_2(0)&=& -4\pi^2\gc-\frac{7\pi^3 i}{3}+2\pi\gc^2 i\\
\end{array}\right.
\end{equation*}
where $\gc=0.5772\dots$ is the Euler constant. We obtain 

\renewcommand{\arraystretch}{2} 
\begin{center}
\begin{tabular}{| p{6.8cm}|}
\hline 
$\begin{array}{ccl}
C_0^{(2,1)}&=& (6\pi-\frac{1}{6}\pi^3-4\pi\gc+\pi\gc^2) \,i\\
C_0^{(3,1)}&=& 2\pi(2-\gc)\,i\\
C_0^{(4,1)}&=& 2\pi i\\
\end{array}$\\
\hline
\end{tabular}\, .
\end{center}

The lattice built on the unique polynomial $q(1/x)=-1/x$ is generated by $q$ itself and the matrix of the exponential torus is $
T_\gl=\diag(1,\gl,\gl,\gl)$. The action of the exponential torus on $I_4+C_0$ results in $T_\gl (I_4+C_0)T_\gl^{-1}=I_4+\gl \, C_0$  and $\ln(I_4+T_\gl C_0T_\gl^{-1})=\gl C_0$. Thus, $\dD_1^+=C_0$ and the alien derivation $\dD_1=C_0$. We can write $\dD_1(\tF)=\tF\, x^L\,C_0\,x^{-L}\,e^{-1/x}$ or $\gD_1(\tF)=\tF\, x^L\,C_0\,x^{-L}$, \ie

\renewcommand{\arraystretch}{2} 
\begin{center}
\begin{tabular}{| p{8.3cm}|}
\hline 
$\begin{array}{ccl}
\gD_1(\tf_2) &=&\displaystyle C_0^{(2,1)}+C_0^{(3,1)}\,\ln x+\frac{1}{2}C_0^{(4,1)}\,\ln^2 x\\
\gD_1(\tf_3) &=& C_0^{(3,1)}+C_0^{(4,1)}\,\ln x\\
\gD_1(\tf_4) &=& C_0^{(4,1)}\,\\
\end{array}$\\
\hline
\end{tabular}\, .
\end{center}

\bibliographystyle{plain}
\nocite{*}
\bibliography{niveau1}
\vfill
\hrule

\bigskip

{\bf Mich\`ele Loday-Richaud}\\
LAREMA - UMR 6093 - 
 2 bd Lavoisier F-49 045 ANGERS cedex 01\\
Laboratoire de Math. - UMR 8628 - 
Universit\'e Paris 11  F-91405 ORSAY\\
{\bf Email}: michele.loday@univ-angers.fr
\bigskip

{\bf Pascal Remy}\\
6 rue Chantal Mauduit F-78 420 Carri\`eres-sur-Seine\\
{\bf Email}: pascal.remy07@orange.fr
\end{document}